\documentclass[12pt, hidelinks]{article}
\usepackage[utf8]{inputenc}
\usepackage{amsmath, amssymb, url, float, graphicx, float, xcolor}
\usepackage{fullpage}
\usepackage[small,bf]{caption}
\usepackage{adjustbox}
\usepackage{hyperref}
\usepackage{subcaption}

\usepackage[style=alphabetic,backend=bibtex]{biblatex}
\bibliography{refs}

\usepackage{booktabs}

\usepackage{tikz, pgfplots}
\pgfplotsset{
    compat=1.16,
}
\usetikzlibrary{topaths,intersections,calc, 3d}
\usetikzlibrary{positioning,chains,fit,shapes,calc,arrows.meta}
\usepgfplotslibrary{fillbetween}

\let\emptyset\varnothing
\let\phi\varphi

\newcommand{\BEAS}{\begin{eqnarray*}}
    \newcommand{\EEAS}{\end{eqnarray*}}
    \newcommand{\BEA}{\begin{eqnarray}}
    \newcommand{\EEA}{\end{eqnarray}}
    \newcommand{\BEQ}{\begin{equation}}
    \newcommand{\EEQ}{\end{equation}}
    \newcommand{\BIT}{\begin{itemize}}
    \newcommand{\EIT}{\end{itemize}}
    \newcommand{\BNUM}{\begin{enumerate}}
    \newcommand{\ENUM}{\end{enumerate}}
    
    \newcommand{\eg}{{\it e.g.}}
    \newcommand{\ie}{{\it i.e.}}

    \newcommand{\ones}{\mathbf 1}
    
    \newcommand{\reals}{\mathbf{R}}
    \newcommand{\integers}{\mathbf{Z}}

    \newcommand{\nullspace}{{\mathcal N}}

    \newcommand{\bmat}[1]{\begin{bmatrix}#1\end{bmatrix}}
    
    \newcommand{\conv}{{\mathop {\bf conv}}}

    \newcommand{\cl}{\mathop{\bf cl}}

    \newcommand{\argmax}{\mathop{\rm argmax}}

    \newif\iftodos

\newcommand{\cone}{\textbf{cone}}

\title{The Convex Geometry of Network Flows}
\author{Theo Diamandis \and Guillermo Angeris}
\date{August 2024}

\begin{document} 
\maketitle 

\begin{abstract}
    In this paper, we derive a number of interesting properties and extensions
    of the convex flow problem from the perspective of convex geometry. We
    show that the sets of allowable flows always can be imbued with a downward
    closure property, which leads to a useful `calculus' of flows, allowing easy
    combination and splitting of edges. We then derive a conic form for the
    convex flow problem, which we show is equivalent to the original problem and
    almost self-dual. Using this conic form, we consider the nonconvex flow
    problem with fixed costs on the edges, \ie, where there is some fixed cost
    to send any nonzero flow over an edge. We show that this problem has almost
    integral solutions by a Shapley--Folkman argument, and we describe a 
    rounding scheme that works well in practice. Additionally, we provide a 
    heuristic for this nonconvex problem which is a simple modification of our
    original algorithm. We conclude by discussing a number of interesting
    avenues for future work.
\end{abstract}

\section{Introduction}
Network flows model systems in a wide range of applications, from flows of cars
to flows of bits. Traditional network flow models assume a linear input-output
relationship for each edge: the flow out of an edge is a linear function of the
flow into that edge, and an extensive literature covers the theory, algorithms,
and applications of these models. (See, \eg,~\cite{ahuja1988network}, 
\cite{williamson2019network}, and references
therein.) Unfortunately, many real-world systems do not exhibit this
linear relationship. For example, in many practical cases, the output of an edge is a concave function of
the input. While nonlinear cost functions have also been
extensively explored in the literature (\eg, see~\cite{bertsekas1998network}
and references therein), nonlinear edge flows---when the flow out of an edge is
a nonlinear function of the flow into it---has received considerably less
attention despite its increased modeling capability.

To model these systems,~\cite{diamandis2024convex} introduced the convex flow
problem, which generalizes the traditional network flow problem to allow for
concave input-output relationships and hypergraph structures. This problem
also generalizes and extends those considered
in~\cite{truemper1978optimal,shigeno2006maximum,vegh2014concave}.
(See~\cite{diamandis2024convex} for additional discussion of this problem and
related literature.) In this paper, we explore the geometry of the convex flow
problem to elucidate a number of interesting properties, many with immediate
practical implications. 

We begin by requiring a downward closure property on the sets of allowable flows.
We prove that requiring this property is equivalent to requiring that the 
utility functions are nondecreasing, as in~\cite{diamandis2024convex}. Having 
downward closed allowable flows has a number of important implications. First, 
downward closure allows for a useful `calculus' of flows, enabling the 
straightforward merging and splitting of network edges. These composition rules 
have immediate practical implications for solvers. 

Second, we derive a conic form for the convex flow problem, which we show is 
equivalent to the original problem. This conic form allows us to derive a number 
of theoretical properties, including a dual that closely resembles the original 
problem, coming close to the self-duality of Bertsekas' extended monotropic 
programming problem~\cite{bertsekas2008extended}.

Finally, the conic problem allows us to extend our analysis to nonconvex flow
problems where each edge has an associated fixed cost for sending any nonzero
flow, a scenario that commonly appears in real-world problems. We show, via
an application of the Shapley--Folkman lemma, that this nonconvex flow problem has
`almost' integral solutions and propose a simple modification to the algorithm
from~\cite{diamandis2024convex} to accommodate these nonconvex scenarios.

\paragraph{Outline.}
In section~\ref{sec:problem}, we define the convex flow problem and discuss
interpretations. We also prove the equivalence between requiring nondecreasing
utility functions and requiring a downward closure condition on the sets of allowable flows. In
section~\ref{sec:calculus}, we show a number of composition rules of downward closed
sets which preserve the downward closure. Then, in section~\ref{sec:conic-form},
we derive an equivalent conic form for the convex flow problem and show that
this conic form is almost self-dual. Finally, we consider the nonconvex flow
problem with fixed costs on the edges in section~\ref{sec:nonconvex}. We show
that this problem is intimately related to the conic form of the convex flow
problem described in the previous section and that it has `almost' integral
solutions.

\section{Problem set up}\label{sec:problem}
In this section, we present the convex flow problem, first defined
in~\cite{diamandis2024convex}, give a few simple properties, and discuss some
important special cases.

\subsection{Problem definition}\label{sec:prob-def}
The \emph{convex flow} problem is the following problem:
\begin{equation}\label{eq:main}
    \begin{aligned}
        & \text{maximize} && {\textstyle U(y) + \sum_{i=1}^m V_i(x_i)}\\
        & \text{subject to} && {\textstyle y = \sum_{i=1}^m A_ix_i}\\
        &&& x_i \in T_i, \quad i=1, \dots, m.
    \end{aligned}
\end{equation}
Here, the \emph{network utility function} $U: \reals^n \to \reals \cup
\{+\infty\}$ and the \emph{edge utility functions} $V_i: \reals^{n_i} \to
\reals \cup \{+\infty\}$ are concave and nondecreasing, the sets $T_i \subseteq
\reals^{n_i}$ are nonempty, closed, and convex, and the matrices $A_i \in \reals^{n \times n_i}$
with $n_i \le n$ are \emph{selector matrices}. Specifically, $A_i$ is a matrix 
of the form
\begin{equation}
    \label{eq:selector}
    A_i = \bmat{a_1 & \dots & a_{n_i}},
\end{equation}
where each $a_k \in \reals^n$ is a distinct unit basis vector. The matrix $A_i$
therefore maps the local indices of edge $i$ to the global indices. We also
assume that an edge need not be used; \ie, that $0 \in T_i$ for all $i$. This
condition makes the proofs simpler and can always be satisfied by appropriately
translating the problem variables and absorbing this translation into the
objective terms.

\begin{figure}[h]
    \centering
    \begin{subfigure}[t]{0.46\textwidth}
        \centering
        \adjustbox{max width=\textwidth}{
            \begin{tikzpicture}
                    \node (v1) at (0,2) {};
                    \node (v2) at (3,2.5) {};
                    \node (v3) at (0,0) {};
                    \node (v4) at (4,-0.5) {};
                
                    \begin{scope}[fill opacity=0.8]
                    \filldraw[fill=blue!70] ($(v1)+(0.5,1)$) 
                        to[out=0,in=180] ($(v3) + (1.5,0)$)
                        to[out=0,in=90] ($(v4) + (1,-0.5)$)
                        to[out=270,in=0] ($0.5*(v3) + 0.5*(v4) + (1,-1)$)
                        to[out=180,in=270] ($(v3) + (-1,0)$)
                        to[out=90,in=180] ($(v1)+(0.5,1)$);
                    \node at (2,-0.5) {\footnotesize $e_1$};
                    \filldraw[fill=red!70] ($(v1)+(0,0.5)$)
                        to[out=0,in=90] ($.5*(v1)+.5*(v3) + (0.55,-0.25)$)
                        to[out=270,in=0] ($(v3) + (0,-1)$)
                        to[out=180,in=270] ($.5*(v1)+.5*(v3) + (-0.55,-0.25)$)
                        to[out=90,in=180] ($(v1)+(0,0.5)$);
                    \node at (0,1) {\footnotesize $e_2$};
                    \filldraw[fill=green!70] ($(v2) + (-.25,1)$)
                        to[out=0,in=90] ($.5*(v2)+.5*(v4) + (0.5,0)$)
                        to[out=270,in=0] ($(v4) + (0.5,-.85)$)
                        to[out=180,in=270] ($.5*(v2)+.5*(v4) + (-0.75,0)$)
                        to[out=90,in=180] ($(v2) + (-.25,1)$);
                    \node at ($.5*(v2)+.5*(v4)$) {\footnotesize $e_3$};
                    \end{scope}
                
                    \foreach \v in {1,2,3,4} {
                        \fill (v\v) circle (0.1);
                    }
                
                    \fill (v1) circle (0.05) node [below] {$v_1$};
                    \fill (v2) circle (0.05) node [above] {$v_2$};
                    \fill (v3) circle (0.05) node [below] {$v_3$};
                    \fill (v4) circle (0.05) node [below] {$v_4$};
                
                \end{tikzpicture}
            }
            \caption{Hypergraph interpretation}
            \label{fig:hypergraph}
        \end{subfigure}
        \hfill
        \begin{subfigure}[t]{0.46\textwidth}
            \centering
            \adjustbox{max width=\textwidth}{
            
                \begin{tikzpicture}[thick,
                    cfmm/.style={draw,circle},
                    every fit/.style={ellipse,draw,inner sep=-2pt,text width=2cm},
                    shorten >= 3pt,shorten <= 3pt,
                ]
                
                \node at (0,2.25) {$S_1$};
                \begin{scope}[yshift=+10mm, start chain=going below, node distance=5mm]
                    \foreach \i in {1,2,3}
                        \node[cfmm,on chain] (e\i) {};
                    \node at ($(e1) + (-0.5,0)$) [left] {\footnotesize $e_1$};
                    \node at ($(e2) + (-0.5,0)$) [left] {\footnotesize $e_2$};
                    \node at ($(e3) + (-0.5,0)$) [left] {\footnotesize $e_3$};
                \end{scope}
                
                \node at (4,2.25) {$S_2$};
                \begin{scope}[xshift=4cm,yshift=+12.5mm,start chain=going below,node distance=5mm]
                    \foreach \i in {1,2,...,4}
                        \node[on chain, circle, draw, fill=black] (v\i) {};
                    \node at ($(v1) + (0.5,0)$) [right] {$v_1$};
                    \node at ($(v2) + (0.5,0)$) [right] {$v_2$};
                    \node at ($(v3) + (0.5,0)$) [right] {$v_3$};
                    \node at ($(v4) + (0.5,0)$) [right] {$v_4$};
                \end{scope}
                
                
                
                \draw[blue] (e1) -- (v1);
                \draw[blue] (e1) -- (v3);
                \draw[blue] (e1) -- (v4);
                
                \draw[red] (e2) -- (v1);
                \draw[red] (e2) -- (v3);
                
                \draw[green] (e3) -- (v2);
                \draw[green] (e3) -- (v4);
        
                \node at (0,-2) {};
                
                \end{tikzpicture}
            }
            \caption{Bipartite graph interpretation}
            \label{fig:bipartite}
        \end{subfigure}

    \caption{A hypergraph with 3 edges and 4 nodes (left) and its
    corresponding bipartite graph representation (right).}
\end{figure}
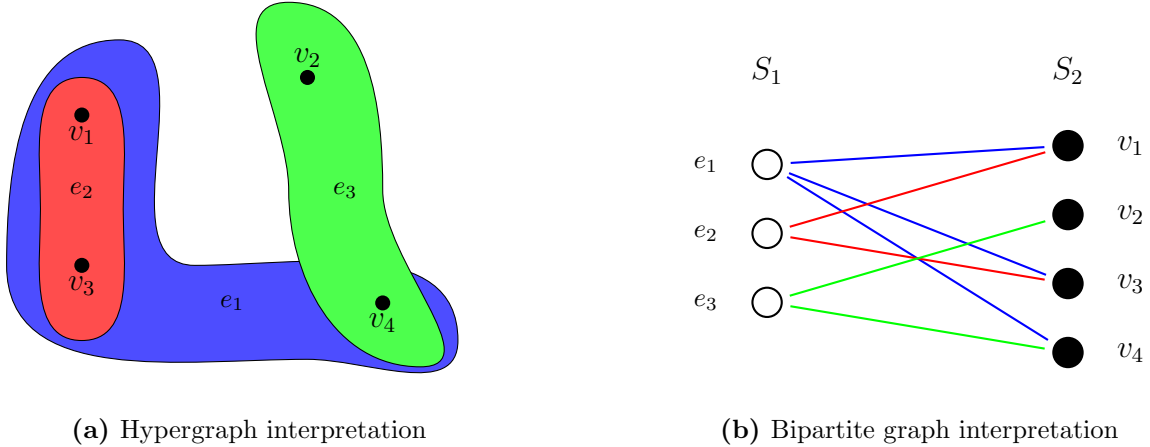

\paragraph{Interpretation: hypergraph flows.}
Following~\cite{diamandis2024convex}, we can interpret this problem as that of
finding the highest-utility allowable flows over a hypergraph with $n$ nodes and
$m$ edges, which may be incident to more than two vertices. The variable $x_i$
indicates the flow over the $i$th hyperedge, which is adjacent to $n_i$
vertices. By convention, we use positive numbers to denote flows out of an edge
(equivalently, into a node) and negative numbers to denote flow into an edge
(equivalently, out of a node). For example, $(x_{i})_k > 0$ is the flow out of
the $i$th edge and into the edge's $k$th incident vertex. Flows $x_i$ over edge
$i$ must lie in $T_i$, the set of allowable flows. The selector matrices map
these edge flows to their corresponding vertices so that $y_j$ is the net flow
at node $j = 1, \dots, n$. Figure~\ref{fig:hypergraph} illustrates this
interpretation.

\paragraph{Interpretation: bipartite graph.}
Alternatively, we may interpret this problem as finding the highest-utility 
flows over a bipartite graph. Here, a set of $m$ vertices $S_1$ is connected to
a set of $n$ vertices $S_2$. We denote the flow from node $i \in S_1$ to its 
$k$th incident node in $S_2$ by $x_{ik}$. The vector $x_i$ is then the set of 
all flows incident to vertex $i \in S_1$, and we it must lie in the set $T_i$.
Again, positive numbers denote flows from $S_1$ to $S_2$, and negative numbers
denote flows in the opposite direction. The selector matrices $A_i$ map the 
flows on edge $i$ to their incident vertices in $S_2$. Thus, $y_j$ is the net 
flow at node $j \in S_2$. Figure~\ref{fig:bipartite} illustrates this 
interpretation.

\paragraph{Examples.}
A number of optimization problems from the literature are special cases
of~\eqref{eq:main}. Classic problems, including maximum flow, minimum cost flow,
and multi-commodity flows, and their generalization to concave edge input-output
relationships, naturally fit into this framework. In engineering, this framework
models the problems of optimal power flow and resource allocation in wireless
networks, which both have nonlinear edge flow relationships. In economics, this
framework includes and generalizes the classic Fisher market model equilibrium
problem. And, in finance, this framework models trading through multiple markets,
where the output of a market is often a concave function of the input---the more
one trades, the worse the price one gets. See~\cite[\S3]{diamandis2024convex}
for details and discussion of these examples.

\subsection{Downward closure and monotonicity}\label{sec:downward-closure}
We say that a set $T \subseteq \reals^n$ is \emph{downward closed} if, for any
$x \in T$ and $x' \le x$, we have $x' \in T$. In other words, if a flow is
feasible, then any smaller flow is also feasible. If $x' \ge x$, we say that the
flow $x'$ \emph{dominates} the flow $x$, since, under any nonnegative utility
function, the flow $x'$ is always at least as `good' as $x$.
In~\cite{diamandis2024convex}, the authors assumed that the functions $U$ and
$\{V_i\}$ in the convex flow problem are nondecreasing. This assumption is, in
fact, equivalent to the sets $\{T_i\}$ being downward closed in the following
sense: if the sets $\{T_i\}$ are downward closed, then the functions $U$ and
$\{V_i\}$ can be replaced with their nondecreasing concave envelopes without
affecting the optimal objective value. Similarly, if the functions $U$ and
$\{V_i\}$ are nondecreasing, then the sets $\{T_i\}$ can be replaced by their
downward closures, \ie,
\[
    \tilde T_i = T_i - \reals^n_+,
\]
without affecting the objective value. This downward closedness property has a
number of immediate and useful implications.

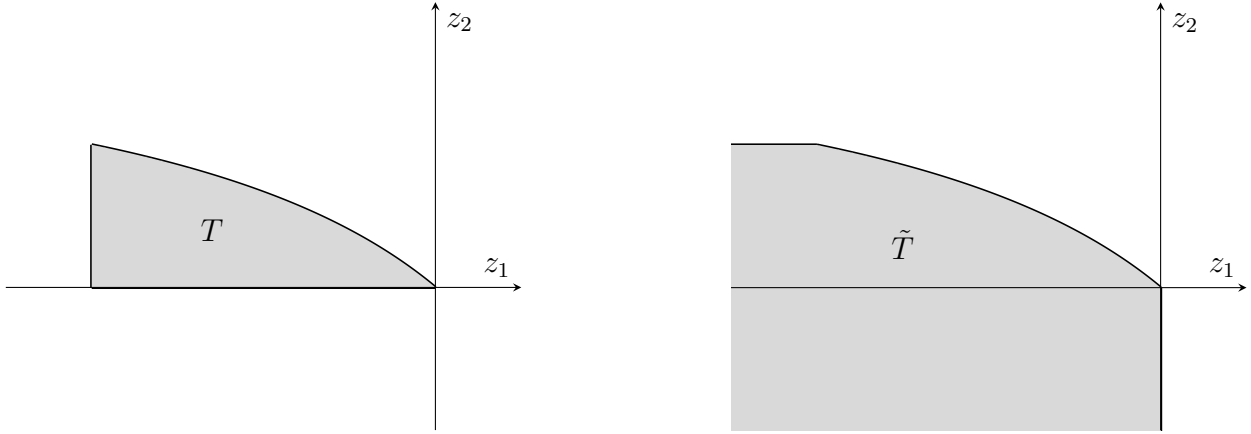
\begin{figure}[h]
    \centering
    \begin{adjustbox}{max width=0.49\textwidth}
        \begin{tikzpicture}[scale=1.0]
    \begin{axis}[
        axis on top=true,
        xmin=-1.25, xmax=0.25,
        ymin=-0.5, ymax=1.0,
        axis lines=center,
        xlabel={$z_1$},
        ylabel={$z_2$},
        grid=major,
        legend pos=outer north east,
        samples=200,
        domain=-5:5,
        xtick=\empty,
        ytick=\empty
        ]                    

        \addplot[black, very thick, domain=-1.00:0.00] {-x / (1 - x)};
        \addplot[black, very thick,] coordinates {(-1, 0.5) (-1, 0)};
        \addplot[black, very thick,] coordinates {(-1, 0) (0, 0)};
        \addplot[gray!30, domain=-1.00:0.00, fill, opacity=0.5, draw=none] {-x / (1 - x)} -| (-1,0) -- cycle;
        
        \node at (-.65, 0.2) {$T$};
    \end{axis}
\end{tikzpicture}
    \end{adjustbox}
    \hfill
    \begin{adjustbox}{max width=0.49\textwidth}
        \begin{tikzpicture}[scale=1.0]
    \begin{axis}[
        axis on top=true,
        xmin=-1.25, xmax=0.25,
        ymin=-0.5, ymax=1.0,
        axis lines=center,
        xlabel={$z_1$},
        ylabel={$z_2$},
        grid=major,
        legend pos=outer north east,
        samples=200,
        domain=-5:5,
        xtick=\empty,
        ytick=\empty
        ]                    

        \addplot[black, very thick, domain=-1.00:0.00] {-x / (1 - x)};
        \addplot[black, very thick,] coordinates {(-1, 0.5) (-2, 0.5)};
        \addplot[black, very thick,] coordinates {(0, 0) (0, -1)};
        \addplot[gray!30, domain=-1.00:0.00, fill, opacity=0.5, draw=none] 
            {-x / (1 - x)} -| (-2,0.5) -| (-1,0) -- cycle;
        \addplot[gray!30, domain=-1.00:0.00, fill, opacity=0.5, draw=none] 
            (0,0) -| (-2,0) -| (-2, -1) -| (0, -1) -- cycle;
        
        \node at (-0.75, 0.15) {$\tilde T$};
    \end{axis}
\end{tikzpicture}
    \end{adjustbox}
    \caption{A set of allowable flows $T$ (left) and its downward closure
    $\tilde T$ (right) for a two-node directed edge that, for input $w$, outputs
    $h(w) = w/(1+w)$ units of flow.}
    \label{fig:downward-closure}
\end{figure}

\paragraph{Example.}
As a simple example, consider a directed edge $i$ with maximum input capacity
$1$ that, when $w$ units of flow enter the edge, outputs $h(w)$ units of
Flow. The corresponding set of allowable flows is
\[
    T_i = \{z \in \reals^2 \mid -1 \le z_1 \le 0 ~~ \text{and} ~~ z_2 \le h(-z_1)\}.
\]
This set is easily verified to be closed and convex, as it is the intersection
of two halfspaces and the hypograph of a concave function, but note that it is
not downward closed.
Figure~\ref{fig:downward-closure} shows a $T_i$ and its downward closure
$\tilde T_i$. The downward closure clearly satisfies the same properties: it is
also closed and convex. We show next that, in a general sense, the set $T$ and
its downward closure $\tilde T$ are `equivalent' for problem~\eqref{eq:main}
if the functions $U$ and $\{V_i\}$ are nondecreasing.

\paragraph{Equivalence proof.}
Consider a finite solution $(y^\star, \{x_i^\star\})$ to~\eqref{eq:main}.
If the objective functions are nondecreasing, then $U(x') \le U(x)$ for any
$x' \le x$, and similarly for the functions $\{V_i\}$. If this solution is not
at the boundary, \ie, if there exists some $x_i^\star$ in the relative interior 
of the corresponding $T_i$, then we can find a nonnegative direction $d \in \reals^{n_i}$
such that $x_i^\star + t d \in T_i$ for some $t > 0$. Since the objective
functions are nondecreasing, this new point will have a objective value equal to
the original solution. As a result, there exists a solution at the boundary, and 
we can replace the sets $\{T_i\}$ with their `downward extension',
\[
    \tilde T_i = T_i - \reals^n_+,
\]
without affecting the solution.

Conversely, if the sets are downward closed and the optimal value is finite,
then, by the downward closure of the $T_i$, there does not exist a nonpositive
direction $d$ such that, for some $t > 0$, $x_i^\star + t d \in T_i$ and the
objective value is larger. (Otherwise, we could find a new point dominated by
$x_i^\star$, \ie, in the downward closure of $T_i$, with a higher objective
value.) Equivalently, all subgradients at this solution must be nonnegative:
$\partial U(y^\star) \subseteq \reals_+^n$ and $\partial V_i(x_i^\star)
\subseteq \reals_+^{n_i}$ for $i= 1, \dots, m$. This fact immediately suggests
that there exists a solution on the boundary and we can replace the objective
functions with their monotonic concave envelopes without affecting the solution.
We will give an alternate proof of this equivalence
in~\S\ref{sec:downward-closure-dual-proof}, after we have derived a dual problem 
for~\eqref{eq:main}.

\section{A calculus of flows}\label{sec:calculus} 
In light of the previous discussion, we will assume that the sets of allowable
flows $\{T_i\}$ are downward closed for the remainder of this paper. We next
discuss a number of properties that directly follow from this condition.
Much of this section generalizes the authors' previous work in the
context of automated market makers~\cite[\S2]{angeris2023geometry}. In the
remainder of this section, we will drop subscripts for convenience.

\paragraph{Definition and interpretation.}
Recall that a set of allowable flows $T$ can be any set satisfying the following 
properties:
\begin{enumerate}
    \item The set $T$ is closed and convex.
    \item The set $T$ is downward closed: if $x \in T$ and $x' \le x$, then $x' \in T$.
    \item The set $T$ contains the zero vector: $0 \in T$.
\end{enumerate}
The three conditions imposed on the set of allowable flows have a natural
interpretation. Convexity means that as more flow enters an edge, the marginal
output does not increase. Downward closure means that positive flow (\ie, flow
out of an edge) can be dissipated. This property often has a nice
interpretation. In power systems, it means that we can dissipate power by, for
example, adding a resistive load. In financial markets, it means that we have
the option to `overpay' for an asset. Finally, the last condition means that we
need not use an edge. This assumption is not fundamental; we can always
translate a set $T$ and absorb the translation into the utility functions. This
assumption, however, will simplify some of the proofs later in this paper.

\subsection{Composition rules}\label{sec:composition}
As a result of the downward closure condition, sets of allowable flows satisfy
certain composition rules. Many of these rules follow directly from the calculus 
of convex sets~\cite[\S2.3]{cvxbook}. For example, the intersection of two sets
of allowable flows yields another set of allowable flows. We discuss a few
important composition rules that will be useful in the rest of this paper below.

\paragraph{Nonnegative matrix multiplication.}
Multiplication of a set of allowable flows by a nonnegative matrix $A \in \reals^{p \times k}$
with $\nullspace(A) = \{0\}$, followed by taking the downward closure, results in
another set of allowable flows:
\[
    AT - \reals^p_+ = \{x \mid x \le Ax' \text{ for some } x' \in T\}.
\]
This resulting set is downward closed by definition, and also closed and convex.
Convexity follows from the fact that convexity is preserved under linear
transforms~\cite[\S2.3.2]{cvxbook} and under downward closure. Closedness of the
set follows from~\cite[Theorem 9.1]{rockafellar1970convex}, as we require $A$ to
be injective. This set has a nice interpretation: given some $x \in T$, each
element of the vector $Ax$ is a weighted `meta-flow' with weights given by the
rows of $A$.

\paragraph{Lifting.}
As a special case of nonnegative matrix multiplication, the lifting of a set
of allowable flows into a larger space is also a set of allowable flows. 
Specifically, let $A$ be a selector matrix (as defined in~\eqref{eq:selector}).
Then the set $AT - \reals^k_+$ is a set of allowable flows in $\reals^{k'}$. This
set describes an edge that connects $k'$ vertices but only allows flow between a
subset of $k$ of them.

\begin{figure}
    \centering
    \begin{adjustbox}{max width=0.3\textwidth}
        \begin{tikzpicture}[scale=1.0]
    \begin{axis}[
        axis on top=true,
        xmin=-1.25, xmax=1.0,
        ymin=-1.25, ymax=1.0,
        axis lines=center,
        xlabel={$z_1$},
        ylabel={$z_2$},
        grid=major,
        legend pos=outer north east,
        samples=200,
        domain=-5:5,
        xtick=\empty,
        ytick=\empty
        ]                    

        \addplot[black, very thick, domain=-1.00:0.00] {-x / (1 - x)};
        \addplot[black, very thick,] coordinates {(-1, 0.5) (-2, 0.5)};
        \addplot[black, very thick,] coordinates {(0, 0) (0, -2)};
        \addplot[gray!30, domain=-1.00:0.00, fill, opacity=0.5, draw=none] 
            {-x / (1 - x)} -| (-2,0.5) -| (-1,0) -- cycle;
        \addplot[gray!30, domain=-1.00:0.00, fill, opacity=0.5, draw=none] 
            (0,0) -| (-2,0) -| (-2, -2) -| (0, -2) -- cycle;
        
        \node at (-0.75, 0.15) {$T$};
    \end{axis}
\end{tikzpicture}
    \end{adjustbox}
    \hfill
    \begin{adjustbox}{max width=0.3\textwidth}
        \begin{tikzpicture}[scale=1.0]
    \begin{axis}[
        axis on top=true,
        xmin=-1.25, xmax=1.0,
        ymin=-1.25, ymax=1.0,
        axis lines=center,
        xlabel={$z_1$},
        ylabel={$z_2$},
        grid=major,
        legend pos=outer north east,
        samples=200,
        domain=-5:5,
        xtick=\empty,
        ytick=\empty
        ]                    

        \addplot[black, very thick, domain=0.00:-1.00] ({-x / (1 - x)}, x);
        \addplot[black, very thick,] coordinates {(0.5, -1) (0.5, -2)};
        \addplot[black, very thick,] coordinates {(-2, 0) (0, 0)};
        \addplot[gray!30, domain=0.00:1.00, fill, opacity=0.5, draw=none] 
            ({x / (1 + x)}, -x) -- (0.5, -2) -- (0,-2) -- cycle;
        \addplot[gray!30, fill, opacity=0.5, draw=none] 
            (0,0) -| (-2,0) -| (-2, -2) -| (0, -2) -- cycle;
        
        \node at (0.2, -.65) {$\tilde T$};
    \end{axis}
\end{tikzpicture}
    \end{adjustbox}
    \hfill
    \begin{adjustbox}{max width=0.3\textwidth}
        \begin{tikzpicture}[scale=1.0]
    \begin{axis}[
        axis on top=true,
        xmin=-1.25, xmax=1.0,
        ymin=-1.25, ymax=1.0,
        axis lines=center,
        xlabel={$z_1$},
        ylabel={$z_2$},
        grid=major,
        legend pos=outer north east,
        samples=200,
        domain=-5:5,
        xtick=\empty,
        ytick=\empty
        ]                    

        \addplot[black, very thick, domain=-1.00:0.00] {-x / (1 - x)};
        \addplot[black, very thick,] coordinates {(-1, 0.5) (-2, 0.5)};
        \addplot[gray!30, domain=-1.00:0.00, fill, opacity=0.5, draw=none] 
            {-x / (1 - x)} -| (-2,0.5) -| (-1,0) -- cycle;
        \addplot[gray!30, domain=-1.00:0.00, fill, opacity=0.5, draw=none] 
            (0,0) -| (-2,0) -| (-2, -2) -| (0, -2) -- cycle;

        \addplot[black, very thick, domain=0.00:-1.00] ({-x / (1 - x)}, x);
        \addplot[black, very thick,] coordinates {(0.5, -1) (0.5, -2)};
        \addplot[gray!30, domain=0.00:1.00, fill, opacity=0.5, draw=none] 
            ({x / (1 + x)}, -x) -- (0.5, -2) -- (0,-2) -- cycle;
        
        \node at (-.5, -.5) {$T + \tilde T$};
    \end{axis}
\end{tikzpicture}
    \end{adjustbox}
    \caption{We take the Minkowski sum of the sets of allowable flows $T$ and
    $\tilde T$ for two directed edges with the form of that in
    figure~\ref{fig:downward-closure} (left and middle) to obtain a new set of
    allowable flows $T + \tilde T$ that corresponds to an undirected edge
    (right).}
    \label{fig:flow-set-addition}
\end{figure}
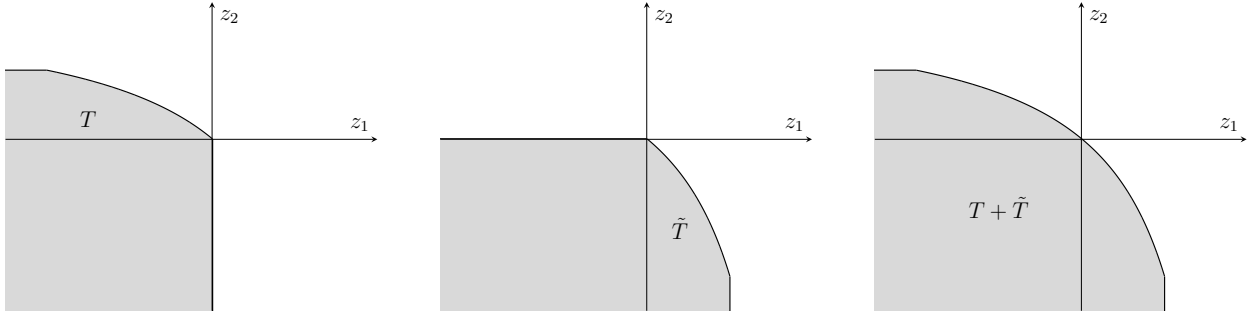

\paragraph{Set addition.}
Finally, under an additional boundedness assumption, the Minkowski sum of
allowable flow sets $T$ and $\tilde T$,
\[
    T + \tilde T = \{x + \tilde x \mid x \in T, \tilde x \in \tilde T\},
\]
is also a set of allowable flows. For this composition rule, we require that
the sets are bounded from above: for a set $T$, there exists some $b$ such that
$x \le b\ones$ for all $x \in T$. This condition means that a bounded input
flow cannot produce infinite output flow. We can interpret the combined set $T
+ \tilde T$ as an aggregate edge that can use either of the two original edges;
see the example in figure~\ref{fig:flow-set-addition}. 

\paragraph{Aggregate edges.}
Using the previous two rules, we can combine edges with possibly non-overlapping 
incident vertices. Importantly, we can view the net flow vector $y$ in~\eqref{eq:main}
as the flow over an `aggregate edge' that connects all vertices with associated
allowable flows
\[
    T = \sum_{i=1}^m \left(A_iT_i - \reals^{n_i}_+\right).
\]
Thus, when the edge utility functions are equal to zero, the convex network
flow problem~\eqref{eq:main} is equivalent to the following problem over one large aggregate
edge:
\[
    \begin{aligned}
        & \text{maximize} && {U(y)}\\
        & \text{subject to} && {\textstyle y \in T}.
    \end{aligned}
\]
While this particular rewriting is not immediately useful, combining or
splitting certain trading sets, for example those with the same incident
vertices, can sometimes help us compute a solution more efficiently. 
(See~\cite[\S4-5]{diamandis2024convex}) for details.)

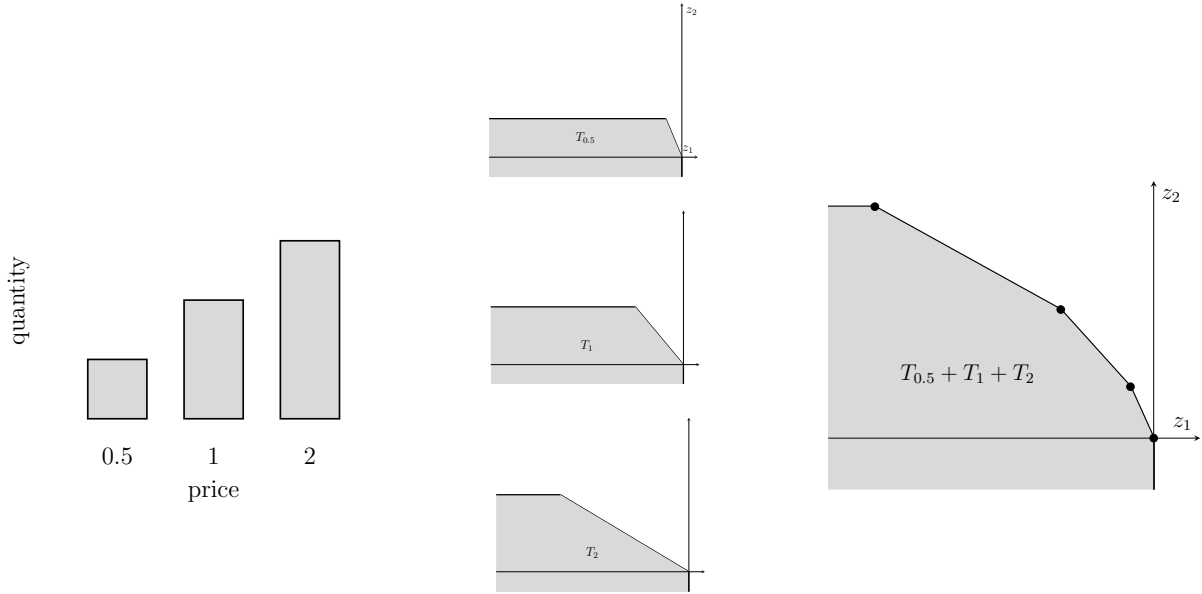
\begin{figure}[h]
    \begin{minipage}{0.3\textwidth}
        \centering
        \vspace{2\baselineskip}
        \begin{adjustbox}{max width=\textwidth}
            \begin{tikzpicture}
\begin{axis}[
    ybar,
    ymin=0,
    ymax=4,
    xmin=0.5,
    xmax=3.5,
    xtick={1,2,3},
    xticklabels={0.5,1,2},
    yticklabels={},
    ylabel=quantity,
    xlabel=price,
    height=6cm,
    ybar=-1cm, 
    enlarge y limits=0.1,
    enlarge x limits=0.5,
    bar width=1cm, 
    axis line style={draw=none}, 
    tick style={draw=none}, 
    enlargelimits=0.05,
]
\addplot[black, fill=gray!30, thick] coordinates {(1,1)};
\addplot[black, fill=gray!30, thick] coordinates {(2,2)};
\addplot[black, fill=gray!30, thick] coordinates {(3,3)};
\end{axis}
\end{tikzpicture}
        \end{adjustbox}        
    \end{minipage}
    \begin{minipage}{0.35\textwidth}
        \centering
        \begin{adjustbox}{max height=0.4\textwidth}
            \begin{tikzpicture}[scale=1.0]
    \begin{axis}[
        axis on top=true,
        xmin=-6, xmax=0.5,
        ymin=-0.5, ymax=4,
        axis lines=center,
        xlabel={$z_1$},
        ylabel={$z_2$},
        grid=major,
        legend pos=outer north east,
        samples=200,
        domain=-4:4,
        xtick=\empty,
        ytick=\empty
        ]
        \addplot[black, very thick, domain=-0.5:0] {-2*x};
        \addplot[gray!30, domain=-0.5:0, fill, opacity=0.6, draw=none] 
            {-2*x} -| (-6, 1) -- cycle;
        \addplot[gray!30, fill, opacity=0.6, draw=none] 
            (0,0) -| (-6, -1) -| cycle;
        
        \node at (-3, 0.5) {$T_{0.5}$};
        \draw[black, very thick] (-0.5, 1) -- (-6, 1);
        \draw[black, very thick] (0, 0) -- (0, -1);

    \end{axis}
\end{tikzpicture}
        \end{adjustbox}
        \vspace{1em}
        \vfill
        \begin{adjustbox}{max height=0.4\textwidth}
            \begin{tikzpicture}[scale=1.0]
    \begin{axis}[
        axis on top=true,
        xmin=-6, xmax=0.5,
        ymin=-0.5, ymax=4,
        axis lines=center,
        xlabel=\empty,
        ylabel=\empty,
        grid=major,
        legend pos=outer north east,
        samples=200,
        domain=-4:4,
        xtick=\empty,
        ytick=\empty
        ]
        \addplot[black, very thick, domain=-1.5:0] {-1*x};
        \addplot[gray!30, domain=-1.5:0, fill, opacity=0.6, draw=none] 
            {-1*x} -| (-6, 1.5) -- cycle;
        \addplot[gray!30, fill, opacity=0.6, draw=none] 
            (0,0) -| (-6, -1) -| cycle;
        
        \node at (-3, 0.5) {$T_{1}$};
        \draw[black, very thick] (-1.5, 1.5) -- (-6, 1.5);
        \draw[black, very thick] (0, 0) -- (0, -1);

    \end{axis}
\end{tikzpicture}
        \end{adjustbox}
        \vspace{1em}
        \vfill
        \begin{adjustbox}{max height=0.4\textwidth}
            \begin{tikzpicture}[scale=1.0]
    \begin{axis}[
        axis on top=true,
        xmin=-6, xmax=0.5,
        ymin=-0.5, ymax=4,
        axis lines=center,
        xlabel=\empty,
        ylabel=\empty,
        grid=major,
        legend pos=outer north east,
        samples=200,
        domain=-4:4,
        xtick=\empty,
        ytick=\empty
        ]
        \addplot[black, very thick, domain=-4:0] {-0.5*x};
        \addplot[gray!30, domain=-4:0, fill, opacity=0.6, draw=none] 
            {-0.5*x} -| (-6, 2) -- cycle;
        \addplot[gray!30, fill, opacity=0.6, draw=none] 
            (0,0) -| (-6, -1) -| cycle;
        
        \node at (-3, 0.5) {$T_{2}$};
        \draw[black, very thick] (-4, 2) -- (-6, 2);
        \draw[black, very thick] (0, 0) -- (0, -1);

    \end{axis}
\end{tikzpicture}
        \end{adjustbox}
    \end{minipage}
    \begin{minipage}{0.3\textwidth}
        \centering
        \vspace{2\baselineskip}
        \begin{adjustbox}{max width=\textwidth}
            \begin{tikzpicture}[scale=1.0]
    \begin{axis}[
        axis on top=true,
        xmin=-7, xmax=1,
        ymin=-1, ymax=5,
        axis lines=center,
        xlabel={$z_1$},
        ylabel={$z_2$},
        grid=major,
        legend pos=outer north east,
        samples=200,
        domain=-4:4,
        xtick=\empty,
        ytick=\empty
        ]
        \addplot[black, very thick, domain=-1:0] 
            (0, -1) -- (0,0) -- (-0.5, 1) -- (-2, 2.5) -- (-6, 4.5)-- (-8, 4.5);
        \addplot[black, mark=*, mark size=2pt] coordinates 
            { (0,0) (-0.5, 1) (-2, 2.5) (-6, 4.5) (-8, 4.5) };
        \addplot[gray!30, draw=none, fill, opacity=0.5] 
            (0, -1) -- (0,0) -- (-0.5, 1) -- (-2, 2.5) -- (-6, 4.5)-- (-8, 4.5) -- (-8, -1) -| cycle;

        
        \node at (-4, 1.25) {$T_{0.5} + T_1 + T_2$};

    \end{axis}
\end{tikzpicture}
        \end{adjustbox}
    \end{minipage}
    \caption{Each tick of the orderbook (left) corresponds to a linear edge with
    a coefficient corresponding to the price (middle). These linear edges
    can be combined into an aggregate edge defining the entire orderbook (right).}
    \label{fig:orderbook}
\end{figure}

\paragraph{Example.}
Often, a directed edge between two nodes has a gain function defined in a 
piecewise manner. For example, consider a financial market between two assets
given by an order book: sellers list the amount of one asset they are willing to
sell for the other at a given price. We can view each `tick' as an individual 
linear edge, which, when combined, define an aggregate edge corresponding to the
entire orderbook. We provide a simple example in figure~\ref{fig:orderbook}.

\section{The conic problem}\label{sec:conic-form}
In this section, we will introduce what looks like a restriction of
problem~\eqref{eq:main}, which we will call the \emph{conic network flow}
problem, defined as
\begin{equation}\label{eq:conic}
    \begin{aligned}
        & \text{maximize} && {\textstyle U(y) + \sum_{i=1}^m V_i(x_i)}\\
        & \text{subject to} && {\textstyle y = \sum_{i=1}^m A_ix_i}\\
        &&& x_i \in K_i, \quad i=1, \dots, m.
    \end{aligned}
\end{equation}
This problem set up is identical to that of~\eqref{eq:main}, except that the
sets $K_i$, instead of being downward closed convex sets,
are downward closed convex cones. A set $K_i$ is called a
\emph{cone} if it satisfies the following property: if $x \in K_i$, then, for
any $\alpha \ge 0$, we must have that $\alpha x \in K_i$. We call any downward
closed convex cone $K_i$ an \emph{allowable flow cone}.

Certainly, every conic flow problem~\eqref{eq:conic} is an instance of a convex
network flow problem~\eqref{eq:main} as every downward closed convex cone is
also, by definition, a downward closed convex set. In this section, we will
show that the converse is also true: every instance of a convex
network flow problem can be turned into an instance of a conic network flow
problem. In this sense, problem~\eqref{eq:main} and problem~\eqref{eq:conic}
are equivalent. We will use the conic problem~\eqref{eq:conic} for the
remainder of this paper to give a number of important theoretical properties,
extensions, and a duality result, all of which easily translate to the
original~\eqref{eq:main}, but are much simpler in the conic formulation.

\subsection{Basic properties}
All of the composition rules presented in~\S\ref{sec:composition} for the
allowable flow sets also hold for the allowable flow cones. More specifically,
given two allowable flow cones (\ie, cones that are downward closed) summation,
intersection, nonnegative scaling, and nonnegative injective matrix
multiplication all yield another allowable flow cone.

\paragraph{Cone is nonpositive.} One immediate consequence of the fact that a
$K\subseteq \reals^d$ is both a cone and downwards closed is that either $K =
\reals^d$ or $K$ contains no strictly positive vectors; that is,
\[
    K \cap \reals_{++}^d = \emptyset.
\]
To see this, let $x \in K$ be any element of $K$ that has only strictly
positive entries $x > 0$. Then for every $d$ vector $y \in \reals^d$, there
exists some $\alpha \ge 0$ such that $y \le \alpha x$. Since $\alpha x$ is in
$K$, as it is a cone, and $K$ is downward closed, then $y \in K$, as required.

\paragraph{Polar cone.} As is standard in convex optimization, given a
cone $K\subseteq \reals^d$ there exists a \emph{polar cone}, defined
\begin{equation}\label{eq:polar-cone}
    K^\circ = \{y \in \reals^d \mid y^Tx \le 0 ~~ \text{for all} ~~ x \in K\}.
\end{equation}
This cone $K^\circ$ is always a closed convex cone (even when $K$ is not). If
$K$ is also a closed convex cone, then we have the following duality result
$(-K^\circ)^\circ = -K$; in other words, the polar cone of the negative polar
cone (called the dual cone) is the negation of the original cone. If, in
addition, the cone $K$ is a downward closed cone with any strictly negative
element (\ie, there is some $x \in K$ with $x < 0$), then we must have that
\[
    K^\circ \cap -\reals_+^d = \{0\}.
\]
(A sufficient condition for this to hold is, \eg, if the cone $K$ has nonempty
interior, which is almost always the case in practice.)

\subsection{Reduction}\label{sec:conic-formulation}
It is clear that the conic problem~\eqref{eq:conic} is a special case of the
original problem~\eqref{eq:main}. In this subsection, we will show that
the converse is also true: any instance of the original problem can be
reduced to an instance of the conic problem.

\paragraph{High level outline.} We begin with an instance of~\eqref{eq:main},
which we write again for convenience:
\begin{equation}\tag{\ref{eq:main}}
    \begin{aligned}
        & \text{maximize} && {\textstyle U(y) + \sum_{i=1}^m V_i(x_i)}\\
        & \text{subject to} && {\textstyle y = \sum_{i=1}^m A_ix_i}\\
        &&& x_i \in T_i, \quad i=1, \dots, m.
    \end{aligned}
\end{equation}
As in~\eqref{eq:main} we have some nondecreasing convex network
utility function $U: \reals^n \to \reals \cup \{-\infty\}$, edge utility
functions $V_i: \reals^{n_i} \to \reals\cup\{-\infty\}$, selector matrices $A_i
\in \reals^{n \times n_i}$, and downward closed sets $T_i \subseteq
\reals^{n_i}$. (From the previous discussion in section~\ref{sec:downward-closure},
we need to require only one of downward closure or montonicity.) 
The variables are the edge flows $x_i \in \reals^{n_i}$ and the
network flows $y \in \reals^n$. Our goal will be to construct some (simple)
nondecreasing net utility function $\tilde U: \reals^{n+1} \to \reals \cup
\{-\infty\}$,
edge utility functions $\tilde V_i: \reals^{n_i+1} \to \reals
\cup \{-\infty\}$, selector matrices $\tilde A_i \in \reals^{(n+1) \times (n_i +
1)}$, and downward closed convex cones $\tilde K_i\subseteq \reals^{n_i+1}$,
such that any solution to the corresponding conic problem~\eqref{eq:conic} over
these new functions, matrices, and sets,
\[
    \begin{aligned}
        & \text{maximize} && {\textstyle \tilde U(\tilde y) + \sum_{i=1}^m \tilde V_i(\tilde x_i)}\\
        & \text{subject to} && {\textstyle \tilde y = \sum_{i=1}^m \tilde A_i\tilde x_i}\\
        &&& \tilde x_i \in \tilde K_i, \quad i=1, \dots, m,
    \end{aligned}
\]
can be (easily) converted to a solution for the original
problem~\eqref{eq:main}. We do this process in two steps. First, we define a
basic cone $\tilde K_i$ associated with each $T_i$, which is essentially the
perspective transformation of $T_i$, done in such a way as to ensure that
$\tilde K_i$ is downward closed when $T_i$ is. We then show that any solution
over this cone, with an additional constraint, always corresponds to a solution
to the original set. Finally, we add this constraint to the objective as an
extra term in the edge cost $\tilde V_i$.

\subsubsection{Flow cone}\label{sec:flow-cone}
We define the flow cone corresponding to a set
$T_i \subseteq \reals^{n_i}$ as
\begin{equation}\label{eq:flow-cone}
    \tilde K_i = \cl\{(x, -\lambda) \in \reals^{n_i}\times \reals \mid x/\lambda \in T_i, ~\lambda > 0\},
\end{equation}
where $\cl$ denotes the closure of a set. This definition is just
the perspective transformation of the set $T_i$, with a sign change in the last
argument. This set is closed (by definition) and convex (see~\cite[\S2.3.3]{cvxbook}). 
It is also downward closed, which we show in appendix~\ref{app:downward-closure-flow-cone}.
We can write the downward closedness of the flow cone in the
following (slightly more useful) way: given any $\lambda'$ such that $-1 \le
\lambda' \le 0$ then
\begin{equation}\label{eq:dom-point}
    (x, \lambda') \in \tilde K_i \quad \text{implies} \quad (x, -1) \in \tilde K_i.
\end{equation}

\paragraph{Recovering the set.} We make use of the following (perhaps obvious)
observation repeatedly. The set $T_i$ can be easily recovered from
the cone $\tilde K_i$ by restricting the last coordinate to be equal to $-1$;
that is,
\begin{equation}\label{eq:recovery}
    T_i = \{x \in \reals^{n_i} \mid (x, -1) \in \tilde K_i\},
\end{equation}
which follows by definition.

\paragraph{Boundary.}
From the flow cone and downward closure, we may define the set $T_i$ via the 
homogenous, nondecreasing, convex function
\[
    \phi_i(x) = \min\{ \lambda \ge 0 \mid (x, -\lambda) \in K_i\}.
\]
Equivalently, we may write $\phi_i$ as the Minkowski functional
\[
    \phi_i(x) = \inf\{\lambda > 0 \mid x/\lambda \in T_i\},
\]
This function's one-level set, $\phi_i(x) = 1$ parameterizes the boundary of $T_i$.
When $\phi_i$ has a nice closed form, this function can be useful in practical applications (see, for
example~\cite{diamandis2023efficient}).

\subsubsection{Rewriting via the flow cone}\label{sec:flow-rewrite}
We start with the original problem~\eqref{eq:main}. Using the cone defined
previously in equation~\eqref{eq:recovery}, we rewrite the original problem
using the flow cone as
\begin{equation}\label{eq:conic-relaxation}
    \begin{aligned}
        & \text{maximize} && {\textstyle U(y) + \sum_{i=1}^m V_i(x_i)}\\
        & \text{subject to} && {\textstyle y = \sum_{i=1}^m A_ix_i}\\
        &&& (x_i, \lambda_i) \in \tilde K_i, \quad i=1, \dots, m\\
        &&& \lambda_i \ge -1,\quad i=1, \dots, m,
    \end{aligned}
\end{equation}
where we have relaxed the constraint that $\lambda_i = -1$ to an inequality.
This relaxation is exact by the `dominating points'
result~\eqref{eq:dom-point}: any solution with $(x_i, \lambda_i)$ may be
replaced with a solution $(x_i, -1)$ which is also feasible, so $x_i \in T_i$
by~\eqref{eq:recovery}. Indeed, in many cases, such as when the set $T_i$ is
locally strictly concave around $0$, one can show that if $\lambda > -1$, there
exists a strictly dominating point $x_i'$ such that $x_i' > x_i$ and $(x_i',
-1) \in \tilde K_i$, so choosing $\lambda_i > -1$ is never optimal. 

\subsubsection{Final rewriting}
Finally, we will take the conic relaxation given
in~\eqref{eq:conic-relaxation}, which, due to the constraint $\lambda_i \ge -1$
is not quite a conic flow problem~\eqref{eq:conic}, and replace the matrices
$A_i$, edge cost functions $V_i$, and variables $(x_i, \lambda_i)$ to get a
problem of the required form.

The first part is easy: let $I(z \ge -1) = 0$ if $z \ge -1$ and $+\infty$ otherwise be
the nonnegative indicator function for a scalar. Note that $I$ is nonincreasing
so $-I$ is nondecreasing, which means we can
rewrite~\eqref{eq:conic-relaxation} by pulling the constraint into the
objective
\begin{equation}\label{eq:conic-constrained-form}
    \begin{aligned}
        & \text{maximize} && {\textstyle U(y) + \sum_{i=1}^m V_i(x_i) - I(\lambda_i \ge -1)}\\
        & \text{subject to} && {\textstyle y = \sum_{i=1}^m A_ix_i}\\
        &&& (x_i, \lambda_i) \in \tilde K_i, \quad i=1, \dots, m.
    \end{aligned}
\end{equation}
Finally, we define the matrix
\[
    \tilde A_i = \bmat{A_i & 0 \\ 0 & 1}.
\]
which is just the matrix $A_i$ with an additional row and column. Setting
$\tilde x_i = (x_i, \lambda_i)$ gives the final result:
\begin{equation}\label{eq:conic-rewriting}
    \begin{aligned}
        & \text{maximize} && {\textstyle \tilde U(\tilde y) + \sum_{i=1}^m \tilde V_i(\tilde x_i)}\\
        & \text{subject to} && {\textstyle \tilde y = \sum_{i=1}^m \tilde A_i\tilde x_i}\\
        &&& \tilde x_i \in \tilde K_i, \quad i=1, \dots, m,
    \end{aligned}
\end{equation}
where we have defined
\[
    \tilde V(x_i, \lambda_i) = V_i(x_i) - I(\lambda_i) 
    \qquad \text{and} \qquad
    \tilde U(y, \tilde \lambda) = U(y).
\]
This problem is exactly of the form of the conic flow
problem~\eqref{eq:conic}, as required.

\subsection{Duality}\label{sec:duality}
Now that we know problem~\eqref{eq:main} and problem~\eqref{eq:conic} are
essentially equivalent (even though the conic problem~\eqref{eq:conic} `seems'
more restrictive) we give a dual reformulation of~\eqref{eq:conic} that is
`almost' self-dual in this section.

\subsubsection{Dual problem}
We will write a simple dual for the conic
problem~\eqref{eq:conic} using standard duality results and a basic
rewriting of the problem.

\paragraph{Lagrangian.} First, we write problem~\eqref{eq:conic}
here again for convenience:
\begin{equation} \tag{\ref{eq:conic}}
    \begin{aligned}
        & \text{maximize} && {\textstyle U(y) + \sum_{i=1}^m V_i(x_i)}\\
        & \text{subject to} && {\textstyle y = \sum_{i=1}^m A_ix_i}\\
        &&& x_i \in K_i, \quad i=1, \dots, m.
    \end{aligned}
\end{equation}
We pull the conic constraint $x_i \in K_i$ into the objective by
defining the indicator functions
\[
    I(x_i \in K_i) = \begin{cases}
        0 & x_i \in K_i\\
        +\infty & \text{otherwise},
    \end{cases}
\]
for $i=1, \dots, m$. We can then rewrite the conic problem as
\[
    \begin{aligned}
        & \text{maximize} && {\textstyle U(y) + \sum_{i=1}^m V_i(x_i) - I(\tilde x_i \in K_i)}\\
        & \text{subject to} && {\textstyle y = \sum_{i=1}^m A_ix_i}\\
        &&& \tilde x_i = x_i, \quad i=1, \dots, m,
    \end{aligned}
\]
where we have introduced the new redundant variables $\tilde x_i \in
\reals^{n_i}$ for each $i=1, \dots, m$. This resulting problem is just a convex
problem with linear constraints. Introducing the Lagrange multipliers $\nu \in
\reals^n$ for the first equality constraint and $\eta_i \in \reals^{n_i}$ for
the second equality constraint, we get the Lagrangian:
\[
    L(x, \tilde x, y, \nu, \eta) = U(y) + \sum_{i=1}^m \left(V_i(x_i) - I(\tilde x_i \in K_i)\right) + \nu^T\left(y - \sum_{i=1}^m A_ix_i\right) + \sum_{i=1}^m \eta_i ^T(x_i - \tilde x_i).
\]

\paragraph{Dual function.} To find the dual function (and therefore the dual
problem) we partially maximize $L$ over the primal variables $x$, $\tilde x$,
and $y$:
\begin{equation}\label{eq:dual-fn}
    g(\nu, \eta) = \bar U(\nu) + \sum_{i=1}^m \bar V_i(\eta_i - A_i^T\nu) + \sum_{i=1}^m \bar I_i(\eta_i).
\end{equation}
Here we have defined
\[
    \bar U(\nu) = \sup_y \left(U(y) + \nu^Ty\right), \quad \bar V_i(\xi_i) = \sup_{x_i} \left(V_i(x_i) + \xi_i^Tx_i\right),
\]
and the functions $\{\bar I_i\}$ as
\[
    \bar I_i(\eta_i) = \sup_{\tilde x_i} \left(-I(\tilde x_i \in K_i) + \tilde x_i^T\eta_i\right),
\]
for each $i=1, \dots, m$. Note that the function $\bar I_i$ is simply the
indicator for the polar cone of $K_i$, defined in~\eqref{eq:polar-cone}. In
other words,
\[
    \bar I_i(\eta_i) = \begin{cases}
        0 & \eta_i \in K_i^\circ\\
        + \infty & \text{otherwise}.
    \end{cases}
\]

\paragraph{Dual problem.} The dual problem is then to minimize the
dual function $g$; \ie,
\[
    \begin{aligned}
        & \text{minimize} && g(\nu, \eta).
    \end{aligned}
\]
When there exists a point in the relative interior of the domain, strong duality holds and, therefore,
the optimal values of the dual problem and the primal problem
are identical. Plugging in the definition of $g$ from~\eqref{eq:dual-fn} into
the objective of the dual problem, and pulling out the indicator
functions $\{\bar I_i\}$ into explicit constraints gives
\[
    \begin{aligned}
        & \text{minimize} && {\textstyle \bar U(\nu) + \sum_{i=1}^m \bar V_i(\eta_i - A_i^T\nu)}\\
        & \text{subject to} && \eta_i \in K_i^\circ, \quad i=1, \dots, m.
    \end{aligned}
\]
We can rewrite the problem to make it more similar to the
original~\eqref{eq:conic}. If we define $\xi_i = A_i^T\nu$ then
\[
    D\nu = \sum_{i=1}^m A_i\xi_i,
\]
where $D$ is a diagonal matrix
\[
    D = \sum_{i=1}^m A_iA_i^T,
\]
with nonnegative diagonal entries. The $j$th diagonal entry, $D_{jj}$, denotes
the degree of node $j$, for $j=1, \dots, n$. The diagonal entries of $D$ are
strictly positive if the hypergraph corresponding to the $A_i$ has no isolated
nodes, or, equivalently, if, for each node $j=1, \dots, n$ there is some edge
$i=1, \dots, m$ such that the $j$th row of $A_i$ is nonzero. In this case,
which we may always assume in practice by removing isolated nodes, the inverse
of $D$ exists so the relationship between $\nu$ and the $\xi_i$ is bijective.
This means we can rewrite the dual problem:
\[
    \begin{aligned}
        & \text{minimize} && {\textstyle \bar U(\nu) + \sum_{i=1}^m \bar V_i(\eta_i - \xi_i)}\\
        & \text{subject to} && {\textstyle D\nu = \sum_{i=1}^m A_i\xi_i}\\
        &&& \eta_i \in K_i^\circ, \quad i=1, \dots, m.
    \end{aligned}
\]
We may absorb the matrix $D$ into the definition of $\bar U$ by replacing $\bar U(\nu)$
with $\bar U(D^{-1}\nu)$ to get the slightly more familiar-looking problem
\begin{equation}\label{eq:conic-dual}
    \begin{aligned}
        & \text{minimize} && {\textstyle \bar U(\nu) + \sum_{i=1}^m \bar V_i(\eta_i - \xi_i)}\\
        & \text{subject to} && {\textstyle \nu = \sum_{i=1}^m A_i\xi_i}\\
        &&& \eta_i \in K_i^\circ, \quad i=1, \dots, m.
    \end{aligned}
\end{equation}
The dual variables may be interpreted as node dual prices $\nu \in \reals^n$ and 
edge dual prices $\xi_i \in \reals^{n_i}$, for $i=1, \dots, n$. We call this
problem the \emph{dual conic flow problem}. Compare this problem~\eqref{eq:conic-dual}
with the original conic flow problem~\eqref{eq:conic}.

\section{Fixed fees}\label{sec:nonconvex}
Finally, we consider the convex network flow problem with fixed fees for the use 
of an edge. In particular, we consider the
following extension of the convex network flow problem~\eqref{eq:main},
which we call the \emph{network flow problem with fees}:
\begin{equation}\label{eq:fixed-fee}
    \begin{aligned}
        & \text{maximize} && {\textstyle U(y) + \sum_{i=1}^m V_i(x_i) + q_i\lambda_i}\\
        & \text{subject to} && {\textstyle y = \sum_{i=1}^m A_ix_i}\\
        &&& (x_i, \lambda_i) \in \{0\} \cup (T_i \times \{-1\}), \quad i=1, \dots, m,
    \end{aligned}
\end{equation}
where the set up and variables are exactly those of the original convex flow
problem~\eqref{eq:main} except with the additional variable $\lambda \in
\reals^m$ and the problem data has the additional fee vector $q \in \reals_+^m$.
This objective is also nondecreasing in all of its variables (as $V_i$ and $U$ 
are, along with the fact that $q \ge 0$).
We note that this problem is not convex since the constraint set is not convex
(in fact, this constraint set is not even connected!) and the problem is NP-hard
to solve, which we prove shortly. However, we will show that this problem can be
approximately solved quite efficiently in practice and is intimately related to
the conic form problem~\eqref{eq:conic} introduced in the previous section.

\paragraph{Interpretation.} The interpretation of the constraint
\begin{equation}\label{eq:hard-constraint}
    (x_i, \lambda_i) \in \{0\} \cup (T_i \times \{-1\})
\end{equation}
is that if $x_i$ is nonzero, then $\lambda_i = -1$. In other words, if we use
edge $i$ by putting any nonzero flow through it, then $x_i \ne 0$ and we are
charged $q_i \ge 0$ for its use. In general, we note that if $q_i > 0$ and $x_i
= 0$, then we will have $\lambda_i = 0$ at optimality, so we may view
$\lambda_i$ as a variable that indicates whether or not edge $i$ is being used.

\paragraph{NP-hardness.}
We show that the network flow problem with fees is NP-hard by reducing the
knapsack problem, which is known to be NP-hard~\cite{karp1972reducibility}, to
an instance of~\eqref{eq:fixed-fee}. The knapsack problem is the following:
given a vector of nonnegative integers $c \in \integers_+^m$ and some integer $b
\ge 0$, find a binary vector $z \in \{0,1\}^m$ such that $c^Tz = b$. This
problem can be reduced to an instance of~\eqref{eq:fixed-fee} with $n = 1$ by
setting $U(y) = y - I(y \ge b)$, $A_i = 1 \in \reals$, $V_i = 0$, $T_i = \{z
\mid z \le c_i\}$, and $q = c$. The problem is
\[
\begin{aligned}
    & \text{maximize} && {-I}(y \ge b) + c^T\lambda\\
    & \text{subject to} && y = {\textstyle \sum_{i=1}^m A_i x_i} \\
    &&& (x_i, \lambda_i) \in \{(0,0)\} \cup ((-\infty,\, c_i] \times \{-1\}), \quad i=1, \dots, m.
\end{aligned}
\]
Note that $c^T(-\lambda) \ge y$ for any feasible point. Since $y$ is
constrained to be at least $b$ then we have that $c^T(-\lambda) \ge b$ for any
feasible point. Finally, maximizing $c^T\lambda$ is the same as minimizing
$c^T(-\lambda) \ge b$, and equality is achieved if and only if there exists
$\lambda \in \{-1, 0\}^m$ such that $c^T(-\lambda) = b$; or, equivalently, when
the optimal objective value of this problem is exactly equal to $-b$. If it
were easy (\ie, polynomial time) to solve this problem, it would be easy to
find a solution to the knapsack problem by verifying that $c^T(-\lambda^\star)
= b$, or to assert that no solution exists if the problem is infeasible or has
optimal objective value larger than $-b$, making this problem at least as hard
as knapsack, which is known to be NP-hard.

\subsection{Integrality constraint}
For the sake of convenience, we will define the set
\[
    Q_i = \{0\} \cup (T_i \times \{-1\}),
\]
such that the constraint~\eqref{eq:hard-constraint} can be written
as
\[
    (x_i, \lambda_i) \in Q_i,
\]
for each $i=1, \dots, m$. In a certain sense, this constraint
encodes the `hard' part of the problem: if the set $Q_i$ were
convex, then the problem would almost be a special case of the original
convex flow problem~\eqref{eq:main}, by pulling the constraint
that $\lambda_i \ge -1$ into the objective.

\paragraph{Convex relaxation.} Given the above discussion, we next examine the 
convex hull of $Q_i$; if we can easily write this convex hull in a compact 
way, then we immediately have a convex relaxation
of the potentially hard problem~\eqref{eq:fixed-fee}. In general, finding the
convex hull of a set may be challenging, \eg, even describing a convex hull can 
require an exponential number of constraints. In this particular special case, we will show that the
convex hull of the set $Q_i$ is intimately related to the flow
cone~\eqref{eq:flow-cone} introduced in the rewriting of the original convex
flow problem~\eqref{eq:main} into the conic flow problem~\eqref{eq:conic}. In
many practical scenarios, finding the flow cone $K_i$ corresponding to the
allowable flows $T_i$ is fairly straightforward, which, in turn, means that
finding the convex hull of $Q_i$ is also fairly straightforward.

\begin{figure}
    \hfill
    \adjustbox{max width=0.45\textwidth}{
        \begin{tikzpicture}
    \begin{axis}[
        axis on top=true,
        view={135}{30}, 
        xmin=-1.4, xmax=0.2,
        ymin=-0.2, ymax=1,
        zmin=-1.5, zmax=0.5, 
        axis lines=center,
        xlabel={$z_1$},
        x label style={anchor=east},
        ylabel={$z_2$},
        zlabel={$\lambda$}, 
        grid=major,
        legend pos=outer north east,
        samples=200,
        domain=-4:4,
        xtick=\empty,
        ytick=\empty,
        ztick=\empty 
        ]
        \addplot3[name path = top, black, variable=t, samples=50, samples y = 0, very thick, domain=0:-1] (t, {sqrt(-t/2)}, -1.0);
        \addplot3[name path = xaxis, draw=none, variable=t, samples=50, samples y = 0, domain=-1:0] (t, 0.0, -1.0);
        \addplot3[gray!30, fill, opacity=0.6, draw=none] 
            fill between[of=top and xaxis];

        \draw[black, very thick] (-1, 0.7071, -1) -- (-1.5, 0.7071, -1);
        \draw[black, very thick] (0, 0, -1) -- (0,-1, -1);
        
        \fill[gray!30, fill, opacity=0.6, draw=none] 
            (0,0,-1) -- (0, -1, -1) -- (-1.5, -1, -1) -- (-1.5, 0, -1) -- (0,0,-1);
        \fill[gray!30, fill, opacity=0.6, draw=none] 
            (-1, 0.7071, -1) -- (-1.5, 0.7071, -1) -- (-1.5, 0, -1) -- (-1, 0, -1) -- cycle;

        \node at (-0.5, -0.1, -1) {$T$};
        \addplot3[black, mark=*, mark size=2] coordinates {(0, 0, 0)};
        \addplot3[black, mark=*, mark size=2] coordinates {(0, 0, -1)} node[below left] {$-1$};
    \end{axis}
\end{tikzpicture}
    }
    \hfill
    \adjustbox{max width=0.45\textwidth}{
        \begin{tikzpicture}
    \begin{axis}[
        axis on top=true,
        view={135}{30}, 
        xmin=-1.5, xmax=0.2,
        ymin=-0.5, ymax=1,
        zmin=-1.2, zmax=0.5, 
        axis lines=center,
        xlabel={$z_1$},
        x label style={anchor=east},
        ylabel={$z_2$},
        zlabel={$\lambda$}, 
        grid=major,
        legend pos=outer north east,
        samples=200,
        domain=-4:4,
        xtick=\empty,
        ytick=\empty,
        ztick=\empty 
        ]
        \addplot3[name path = top, black, variable=t, samples=50, samples y = 0, very thick, domain=0:-1] (t, {sqrt(-t/2)}, -1.0);
        \addplot3[name path = mid, black, variable=t, samples=50, samples y = 0, domain=0:-1] (t, {0.5*sqrt(-t)}, -0.5);
        \addplot3[name path = xaxis, draw=none, variable=t, samples=50, samples y = 0, domain=-1:0] (t, 0.0, 0.0);
        \addplot3[gray!30, fill, opacity=0.6, draw=none] 
            fill between[of=mid and xaxis];
        \addplot3[gray!30, fill, opacity=0.6, draw=none] 
            fill between[of=mid and top];

        \draw[black, very thick] (-1, 0.7071, -1) -- (-1.5, 0.7071, -1);
        \draw[black] (-1, 0.5, -0.5) -- (-1.5, 0.5, -0.5);
        \draw[black, very thick] (0, 0, -1) -- (0,-0.5, -1);
        \draw[black] (0, 0, -0.5) -- (0,-0.5, -0.5);
        
        \fill[gray!30, fill, opacity=0.6, draw=none] 
            (0,0,0) -- (0, -0.5, 0) -- (-1.5, -0.5, 0) -- (-1.5, 0, 0) -- (0,0,0);
        \fill[gray!30, fill, opacity=0.6, draw=none] 
            (-1, 0.7071, -1) -- (-1.5, 0.7071, -1) -- (-1.5, 0.5, -0.5) -- (-1, 0.5, -0.5) -- cycle;
        \fill[gray!30, fill, opacity=0.6, draw=none] 
            (-1, 0.5, -0.5) -- (-1.5, 0.5, -0.5) -- (-1.5, 0, 0) -- (-1, 0, 0) -- cycle;
        
        \fill[gray!30, fill, opacity=0.6, draw=none] 
            (0, 0, 0) -- (0,0,-1) -- (0, -0.5, -1) -- (0, -0.5, 0) -- cycle;
        
        \draw[black, dashed] (0,0,0) -- (-0.125, 0.25, -1);
        \draw[black, dashed] (0,0,0) -- (-0.5, 0.5, -1);
        \draw[black, dashed] (0,0,0) -- (-1, 0.7071, -1);
        \draw[black, dashed] (0,0,0) -- (-1.4, .7, -.7);
        \draw[black, dashed] (0,0,0) -- (-1.5, .5, -.5);
        \draw[black, dashed] (0, 0, 0) -- (0, -0.3, -1);

        \node at (-0.5, -0.1, -1) {$T$};
        \addplot3[black, mark=*, mark size=2] coordinates {(0, 0, 0)};
    \end{axis}
\end{tikzpicture}
    }
    \hfill
    \caption{The set $Q_i$ (left) and its convex hull $\conv(Q_i)$ (right).
    }
    \label{fig:edge-convex-hull}
\end{figure}
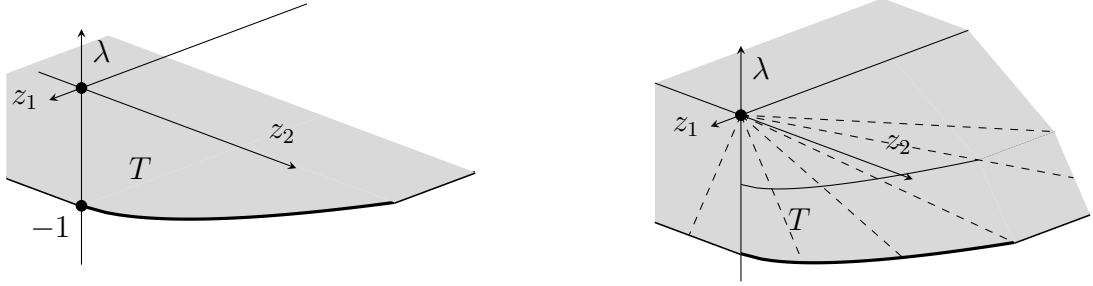

\subsection{Convex hull}
Given the above discussion, we will show that the convex hull of $Q_i$, written
$\conv(Q_i)$, is equal to all elements of the corresponding flow cone $K_i$
when the elements' last entry (corresponding to $\lambda_i$) lies between $-1$
and $0$. Written out we will show that:
\[
    \conv(Q_i) = K_i \cap (\reals^n \times [-1, 0]).
\]
See figure~\ref{fig:edge-convex-hull} for an example.
As a reminder, the flow cone $K_i \subseteq \reals^{n+1}$ for a given allowable
flow set $T_i \subseteq \reals^n$ is defined, using~\eqref{eq:flow-cone}, as
\[
    K_i = \cl\{(x, -\lambda) \in \reals^n \times \reals \mid x/\lambda \in T_i, ~ \lambda > 0\}.
\]
To simplify notation, we define
\begin{equation}\label{eq:clipped-cone}
    \bar K_i = K_i \cap (\reals^n \times [-1, 0]),
\end{equation}
which is the cone $K_i$ with the last element restricted to lie between $-1$
and $0$. Of course, this set is also convex as it is the intersection of two
convex sets.

\paragraph{Reverse inclusion.} First, we show the reverse inclusion: that $\bar
K_i \subseteq \conv(Q_i)$. Let $(x, -\lambda) \in \bar K_i$ (note the negative
here) with $\lambda > 0$, then, we will show that $(x, -\lambda)$ can be
written as the convex combination of an element in $T_i \times \{-1\}$ and the
zero vector,
and so also must lie in $\conv(Q_i)$. By definition,
if $(x, -\lambda) \in \bar K_i$, then $x/\lambda \in T_i$ for $0 < \lambda \le
1$. But, this is the same as saying
\[
    (x/\lambda, -1) \in Q_i.
\]
Finally, since $0 \in Q_i$, then any convex combination of $0$ and $(x/\lambda,
-1)$ is in the convex hull of $Q_i$. Thus,
\[
    (x, -\lambda) = \lambda(x/\lambda, -1) + (1-\lambda)0 \in \conv(Q_i),
\]
so long as $\lambda > 0$. On the other hand, if $\lambda = 0$, then
we know that $x/\lambda' \in T_i$ for all $\lambda' > 0$, so
\[
    (x, \lambda') = \lambda'(x/\lambda', -1) + (1-\lambda')0 \in \conv(Q_i).
\]
Sending $\lambda' \to 0$ gives the result, since $Q_i$ is closed as it is the
union of two closed sets. Putting it all together, this implies that
$\conv(Q_i) \supseteq \bar K_i$.

\paragraph{Forward inclusion.} Now, we show the forward inclusion: that
$\conv(Q_i) \subseteq \bar K_i$. Note that $Q_i \subseteq \bar K_i$ since, by
definition
\[
    T_i \times \{-1\} \subseteq \bar K_i,
\]
and $0 \in \bar K_i$, also essentially by
definition. This immediately implies that
\[
    \conv(Q_i) \subseteq \conv(\bar K_i) = \bar K_i
\]
where we have used the fact that the convex hull of a convex set is itself.

\paragraph{Discussion.} Putting the above two points together, we get the claim
that
\begin{equation}\label{eq:fixed-cost-equiv}
    \conv(Q_i) = \bar K_i.
\end{equation}
In other words, the convex hull of the `hard' set is exactly the flow cone $K_i$
with the additional constraint that the last entry must be restricted to lie
between $0$ and $-1$. One interesting interpretation of this claim is that we
may view the cone $K_i$ as the \emph{conic completion} of the set $Q_i$. More
generally, $\cone(Q_i)$ is defined as the set containing all conic (\ie,
nonnegative) combinations of the elements of $Q_i$. Since is it not hard to
show that $\cone(Q_i) = \cone(\conv(Q_i))$, we have
\[
    \cone(Q_i) = \cone(\conv(Q_i)) = \cone(\bar K_i) = K_i,
\]
where the second equality follows from~\eqref{eq:fixed-cost-equiv}, while the
last simply follows from definitions.

\subsection{Convex relaxation}
Using the result~\eqref{eq:fixed-cost-equiv} derived in the previous section, a convex relaxation
of the network problem with fees~\eqref{eq:fixed-fee} is
\begin{equation}\label{eq:fixed-fee-relaxation}
    \begin{aligned}
        & \text{maximize} && {\textstyle U(y) + \sum_{i=1}^m \left(V_i(x_i) + q_i\lambda_i\right)}\\
        & \text{subject to} && {\textstyle y = \sum_{i=1}^m A_ix_i}\\
        &&& (x_i, \lambda_i) \in K_i, ~ {-1} \le \lambda_i \le 0, \quad i=1, \dots, m.
    \end{aligned}
\end{equation}
Note that we have replaced the nonconvex constraint $(x_i,
\lambda_i) \in Q_i$ with the convex constraint $(x_i, \lambda_i) \in
\conv(Q_i)$, or, equivalently, using the facts derived in the previous section,
replaced it with the constraint that $(x_i, \lambda_i) \in K_i$ and $-1 \le
\lambda_i \le 0$.

\paragraph{Conic formulation.} This convex relaxation is also a special case of
the conic flow problem~\eqref{eq:conic} in a very natural way. First, note that
the constraint that $\lambda \le 0$ is redundant using the definition of $K_i$.
We may then pull the remaining constraint on $\lambda_i$, that $\lambda_i \ge -1$
into an indicator function, $I(\lambda_i \ge -1)$ and place it in the objective.
This indicator function is nonincreasing, so its negation is nondecreasing,
and we get the final problem
\begin{equation}\label{eq:conic-fees-form}
    \begin{aligned}
        & \text{maximize} && {\textstyle U(y) + \sum_{i=1}^m \left(V_i(x_i) + q_i\lambda_i - I(\lambda_i \ge -1)\right)}\\
        & \text{subject to} && {\textstyle y = \sum_{i=1}^m A_ix_i}\\
        &&& (x_i, \lambda_i) \in K_i, \quad i=1, \dots, m.
    \end{aligned}
\end{equation}
If we use the same trick used in~\S\ref{sec:conic-formulation} to rewrite
the matrices $A_i$, we receive an instance of the conic flow
problem~\eqref{eq:conic} since the objective is nondecreasing and the $\{K_i\}$
are cones. Indeed, the formulation found here~\eqref{eq:conic-fees-form} is
essentially identical to the formulation found
in~\eqref{eq:conic-constrained-form}, except with the addition of the fixed
costs $q \ge 0$.

\paragraph{Integrality gap.} In the previous argument in~\S\ref{sec:flow-rewrite},
we used the fact that, if $\lambda_i > -1$, then we could set $\lambda_i = -1$
and always remain feasible with no change in objective value. However in 
problem~\eqref{eq:conic-fees-form}, the objective value
would decrease, so this argument does not apply. We will show next that
we expect the solution to be close to integral in the special case that the $V_i = 0$, which
is common in practice (see, for example, the applications in~\cite[\S3]{diamandis2024convex}).

\subsection{Tightness of the relaxation}
Just how tight do we expect the relaxation
to be? We will show that in the case that $V_i = 0$, if $m$, the number of
edges, is much larger than $n$, the number of nodes, then most of the
$\lambda_i$ will be integral. More specifically, we will show that, given any
solution to the relaxation, we can recover a solution such that at least
$m-n-1$ indices $i$ satisfy $(x_i, \lambda_i) \in Q_i$. If $m \gg n$, \ie, the
number of nodes is much smaller than the number of edges, as is usually the
case in practice, then most of the solution is integral.

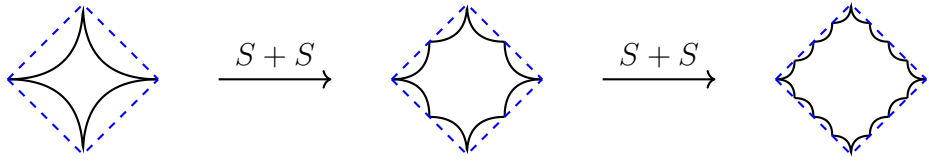
\begin{figure}[t]
    \centering
    \adjustbox{max width=0.85\textwidth}{
        \hfill
\begin{tikzpicture}[scale=1.0,baseline=-0.5ex]
    \draw[black, thick, domain=-1:1, samples=300] plot (\x, {(1 - sqrt(abs(\x)))^2});
    \draw[black, thick, domain=-1:1, samples=300] plot (\x, {-(1 - sqrt(abs(\x)))^2});
    \draw[blue, dashed, thick, domain=-1:1, samples=100] plot (\x, {1 - (abs(\x))});
    \draw[blue, dashed, thick, domain=-1:1, samples=100] plot (\x, {(abs(\x)) - 1});
\end{tikzpicture}
\hspace{0.5cm}
\hfill
\begin{tikzpicture}[baseline=-0.5ex]
    \draw[->, black, thick] (0,0) -- (1.5,0) node[midway, above] {$S + S$};
\end{tikzpicture}
\hfill
\hspace{0.5cm}
\begin{tikzpicture}[scale=0.5, baseline=-0.5ex]
    \draw[black, thick, domain=-1:1, samples=200] plot (\x, {1 + (1 - sqrt(abs(\x)))^2});
    \draw[black, thick, domain=0:1, samples=200] plot (1 + \x, {(1 - sqrt(abs(\x)))^2});
    \draw[black, thick, domain=0:-1, samples=200] plot (-1 + \x, {(1 - sqrt(abs(\x)))^2});

    \draw[black, thick, domain=-1:1, samples=200] plot (\x, {-1 - (1 - sqrt(abs(\x)))^2});
    \draw[black, thick, domain=0:1, samples=200] plot (1 + \x, {-(1 - sqrt(abs(\x)))^2});
    \draw[black, thick, domain=0:-1, samples=200] plot (-1 + \x, {-(1 - sqrt(abs(\x)))^2});

    \draw[blue, dashed, thick, domain=-2:2, samples=100] plot (\x, {2 - (abs(\x))});
    \draw[blue, dashed, thick, domain=-2:2, samples=100] plot (\x, {(abs(\x)) - 2});
\end{tikzpicture}
\hfill
\hspace{0.5cm}
\begin{tikzpicture}[baseline=-0.5ex]
    \draw[->, black, thick] (0,0) -- (1.5,0) node[midway, above] {$S + S$};
\end{tikzpicture}
\hfill
\hspace{0.5cm}
\begin{tikzpicture}[scale=0.25, baseline=-0.5ex]
    \draw[black, thick, domain=-1:1, samples=200] plot (\x, {3 + (1 - sqrt(abs(\x)))^2});
    \draw[black, thick, domain=0:1, samples=200] plot (1 + \x, {2 + (1 - sqrt(abs(\x)))^2});
    \draw[black, thick, domain=0:-1, samples=200] plot (-1 + \x, {2 + (1 - sqrt(abs(\x)))^2});
    \draw[black, thick, domain=0:1, samples=200] plot (2 + \x, {1 + (1 - sqrt(abs(\x)))^2});
    \draw[black, thick, domain=0:-1, samples=200] plot (-2 + \x, {1 + (1 - sqrt(abs(\x)))^2});
    \draw[black, thick, domain=0:1, samples=200] plot (3 + \x, {(1 - sqrt(abs(\x)))^2});
    \draw[black, thick, domain=0:-1, samples=200] plot (-3 + \x, {(1 - sqrt(abs(\x)))^2});

    \draw[black, thick, domain=-1:1, samples=200] plot (\x, {-(3 + (1 - sqrt(abs(\x)))^2)});
    \draw[black, thick, domain=0:1, samples=200] plot (1 + \x, {-(2 + (1 - sqrt(abs(\x)))^2)});
    \draw[black, thick, domain=0:-1, samples=200] plot (-1 + \x, {-(2 + (1 - sqrt(abs(\x)))^2)});
    \draw[black, thick, domain=0:1, samples=200] plot (2 + \x, {-(1 + (1 - sqrt(abs(\x)))^2)});
    \draw[black, thick, domain=0:-1, samples=200] plot (-2 + \x, {-(1 + (1 - sqrt(abs(\x)))^2)});
    \draw[black, thick, domain=0:1, samples=200] plot (3 + \x, {-((1 - sqrt(abs(\x)))^2)});
    \draw[black, thick, domain=0:-1, samples=200] plot (-3 + \x, {-((1 - sqrt(abs(\x)))^2)});

    \draw[blue, dashed, thick, domain=-4:4, samples=100] plot (\x, {4 - (abs(\x))});
    \draw[blue, dashed, thick, domain=-4:4, samples=100] plot (\x, {(abs(\x)) - 4});
\end{tikzpicture}
\hfill
    }
    \caption{
        A visual representation of the Shapley--Folkman lemma for the (nonconvex) $1/2$-norm
        ball. As we take the Minkowski sum of the set with itself, it becomes
        closer and closer to its convex hull, the $1$-norm ball.
    }
    \label{fig:shapley-folkman}
\end{figure}

\paragraph{Shapley--Folkman lemma.} We state the Shapley--Folkman lemma here
in its standard form. Let $S_1, \dots, S_m \subseteq \reals^{n+1}$ be any subsets
(convex or nonconvex) of $\reals^{n+1}$. Then, for any
\[
    y = x_1 + \dots + x_m,
\]
where $x_i \in \conv(S_i)$ for $i=1, \dots, m$, there exists $\tilde x_i \in
\conv(S_i)$ with $i=1, \dots, m$, such that
\[
    y = \tilde x_1 + \dots + \tilde x_m,
\]
which satisfy $\tilde x_i \in S_i$ for at least $m-n-1$ indices $i$. In other
words, given any vector $y$, which lies in the (Minkowski) sum of the convex
hulls of the $S_i$, we can find $\tilde x_i$, which sum to $y$, such that at
least $m-n-1$ lie in the original sets $S_i$, while the remainder lie in the
convex hull, $\conv(S_i)$. Intuitively, this lemma states that the sum of convex
sets becomes closer and closer to its convex hull as the number of sets gets
large. See figure~\ref{fig:shapley-folkman} for an example.

Given a solution to the convex
relaxation~\eqref{eq:conic-fees-form}, we will use this lemma to construct a
solution that has the same objective value as the original, yet almost all
penalties $\lambda_i$ will be integral: either $-1$ or $0$. This will then let
us construct an approximate solution to~\eqref{eq:fixed-fee} and bound the
difference between the optimal objective value of~\eqref{eq:fixed-fee} and
that of the relaxation~\eqref{eq:fixed-fee-relaxation}.

\paragraph{Constructing an approximate solution.} Assume we are given feasible
flows and penalties for the relaxation; \ie, we have a solution $\{(x_i^\star,
\lambda_i^\star)\}$ to the relaxation of the fixed-fee
problem~\eqref{eq:fixed-fee-relaxation}. From this solution to the
relaxation~\eqref{eq:fixed-fee-relaxation}, we will construct a feasible point
for the original fixed-fee problem~\eqref{eq:fixed-fee} which we will then show
is `close' to the optimal value, under certain conditions. We write the exact
problem we are considering (the special case of~\eqref{eq:fixed-fee-relaxation}
when $V_i = 0$) for convenience:
\[
    \begin{aligned}
        & \text{maximize} && {\textstyle U(y) + q^T\lambda}\\
        & \text{subject to} && {\textstyle y = \sum_{i=1}^m A_ix_i}\\
        &&& (x_i, \lambda_i) \in \bar K_i, \quad i=1, \dots, m.
    \end{aligned}
\]
Here, we have used the definition of $\bar K_i$ from~\eqref{eq:clipped-cone},
and, from the previous discussion~\eqref{eq:fixed-cost-equiv}, we know that
$\bar K_i = \conv(Q_i)$. 
First, note that, by definition
\[
    y^\star = \sum_{i=1}^m A_i x_i^\star,
\]
and that $(x_i^\star, \lambda_i^\star) \in \conv(Q_i)$ for $i=1, \dots, m$.
Now, from the problem statement, we have
\[
    c = q^T\lambda^\star,
\]
with $q \ge 0$, where $c \le 0$ stands for the `fixed cost' part of the
objective. Rewriting slightly,
\[
    \bmat{y^\star \\ c} =
    \sum_{i=1}^m \bmat{
            A_i & 0\\
            0 & q_i
    } \bmat{x_i^\star \\ \lambda_i^\star},
\]
where $(x_i^\star, \lambda_i^\star) \in \bar K_i$, or, equivalently,
using~\eqref{eq:fixed-cost-equiv}, $(x_i^\star, \lambda_i^\star) \in
\conv(Q_i),$ for $i=1, \dots, m$. From the Shapley--Folkman lemma, there exist
$(\tilde x_i^\star, \tilde\lambda_i^\star) \in \conv(Q_i)$ such that
\[
    \bmat{y^\star \\ c} =
    \sum_{i=1}^m \bmat{
            A_i & 0\\
            0 & q_i
    } \bmat{\tilde x_i^\star \\ \tilde \lambda_i^\star},
\]
and at least $m-n - 1$ indices $i$ satisfy $(\tilde x_i^\star, \tilde
\lambda_i^\star) \in Q_i$. In other words, if $m \gg n$, we have an `almost'
feasible solution for the original problem~\eqref{eq:fixed-fee}, except for $n+1$
indices $i$. Consider these nonintegral indices, of which there are at most $n+1$. In
this case, we know, from the dominated point condition~\eqref{eq:dom-point} for
$K_i$ (and therefore for $\bar K_i$) that if $(\tilde x_i^\star, \tilde
\lambda_i^\star) \in \bar K_i$, then $(\tilde x_i^\star, -1) \in \bar K_i$. But
we know that $(\tilde x_i^\star, -1) \in Q_i$, making this point also feasible
for the original problem~\eqref{eq:fixed-fee}. In English: if we were charged
less than the full amount due to the relaxation (\ie, we were charged
$q_i\lambda_i$), we can always choose to be charged the full amount for the
same flow ($-q_i$) and be feasible for the original problem. This means that a
feasible solution for the original problem~\eqref{eq:fixed-fee} will be to set
\begin{equation}\label{eq:feasible-approx}
    (x_i^0, \lambda_i^0)  = \begin{cases}
        (\tilde x_i^\star, \tilde \lambda_i^\star) & (\tilde x_i^\star, \tilde \lambda_i^\star) \in Q_i\\
        (\tilde x_i^\star, -1) & \text{otherwise},
    \end{cases}
\end{equation}
for $i=1, \dots, m$. Note that $(x_i^0, \lambda_i^0) \in Q_i$ for each $i$ and
so is a feasible point for the original fixed-fee problem~\eqref{eq:fixed-fee},
leading to the same net flows $y^\star$ as the solution to the relaxation, but
the cost incurred, $q^T\lambda^0$ differs by at most
\[
    q^T(\lambda^0 - \tilde\lambda^\star),
\]
from that of the relaxation, $q^T\tilde \lambda^\star = c$. Since most entries
of $\lambda^0 - \tilde \lambda^\star$ are zero, by the previous argument, then
we expect this cost to be small. We give a simple bound on this, and therefore
in the objective gap between the relaxation and the original problem, in what
follows.

\paragraph{Bounding the optimal objective value.} Let $p^0$ be the optimal
objective value for the relaxation~\eqref{eq:fixed-fee-relaxation} and let
$p^\star$ be the optimal objective value for the original
problem~\eqref{eq:fixed-fee}. Then, since~\eqref{eq:fixed-fee-relaxation} is a
relaxation of~\eqref{eq:fixed-fee}, we know that
\[
    p^\star \le p^0.
\]
From the previous discussion, we have a feasible
point~\eqref{eq:feasible-approx} for the original problem. By construction, we
know that the net flows $y^\star$ remain unchanged, so the net utility
$U(y^\star)$ in the objective similarly remains unchanged. On the other hand,
the cost incurred $q^T\lambda^0$ is larger than that of the relaxation, $c$, by
$q^T(\tilde \lambda^\star -\lambda^0)$, so we have the following bound
\[
    p^0 + q^T(\lambda^0 - \tilde\lambda^\star) \le p^\star \le p^0.
\]
If we solved the relaxation, then we immediately have a two-sided bound on the
optimal objective value as given above. On the other hand, we can give a simple
bound that does not require solving the relaxation. Since we know that
at most $n+1$ entries of $\lambda^0$ will differ from those of $\tilde\lambda^\star$
by Shapley--Folkman, then
\[
    p^0 - (n+1)\left(\max_i q_i\right) \le p^\star \le p^0.
\]
Or, equivalently
\[
    0 \le p^0 - p^\star \le (n+1)(\max_i q_i).
\]

\subsection{Fixed cost dual problem}\label{sec:fixed-cost-dual}
Finally, we derive the dual problem of the fixed-fee problem~\eqref{eq:fixed-fee}
with zero edge costs
and show that the algorithm developed in~\cite{diamandis2024convex} can still be
applied to this problem with minimal modifications. In fact, we lose little
computational efficiently by solving the fixed-fee problem directly; of course,
the solution is not guaranteed to be optimal.

\paragraph{Dual function.}
Using a similar derivation to that in~\S\ref{sec:duality}, we write the
dual function as
\[
    g(\nu) 
    = \bar U(\nu) + 
    \sum_{i=1}^m \sup_{(x_i, \lambda_i) \in Q_i} \left((A_i^T\nu)^Tx_i + \lambda_i q_i\right).
\]
Note that the support function over $Q_i$, \ie, the expression in the sum, can
be easily evaluated: we simply compute the support function of $T_i$ less $q_i$
and compare this value to $0$. If we define
\begin{equation}\label{eq:support-simple}
    f_i(\xi) = \sup_{x_i \in T_i} \xi^Tx_i,
\end{equation}
then we can write
\[
    \sup_{(x_i, \lambda_i) \in Q_i} \left(\xi^Tx_i + \lambda_i q_i\right)
    = \max\left\{f_i(\xi) - q_i,\, 0 \right\}.
\]
The optimal points for the support function are as follows. If $f_i(\xi) \ge
q_i$, then $\lambda_i^\star = -1$ and $x_i^\star$ can be any solution to the
subproblem~\eqref{eq:support-simple}. On the other hand, if $f_i(\xi) \le q_i$,
then $\lambda_i^\star = 0$ and $x_i^\star = 0$ is a solution. (There may be
many solutions if $f_i(\xi) = q_i$.) This observation allows us to apply the
algorithm from~\cite{diamandis2024convex} `off the shelf' to solve the dual
problem. 

\paragraph{Primal feasibility.}
By construction, using the $\lambda_i^\star$ and $x_i^\star$ from
above always results in integral solutions to the
edge subproblem. Of course, these solutions may not be primal feasible: the
net flows $y^\star$ solving the $\bar U$ subproblem may not be equal to the sum of the
solutions $x_i^\star$ to the edge subproblems; \ie, it is possible
that
\[
    y^\star = \sum_{i=1}^m A_ix_i^\star,
\]
is not true. However, we can always construct a primal feasible net flow $\hat
y$ using these subproblem solutions $\{x_i^\star\}$ by setting
\[
    \hat y = \sum_{i=1}^m A_i x_i^\star.
\]

\paragraph{Verifying optimality.}
We can verify the optimality of the primal feasible point $(\hat y, \{x_i^\star\})$
by checking it's objective value. From duality, we know that
\[
    g(\nu^\star) = U(y^\star) \ge U(\hat y).
\]
If these values are equal, we know that the point $(\hat y, \{x_i^\star\})$ is
also optimal.

\subsection{Numerical experiments}
Here, we present some simple numerical experiments to illustrate the results of
this section. We emphasize the these experiments are by no means exhaustive, but
we leave a more thorough investigation to future work.

\paragraph{Setup.}
We consider the order routing problem from~\cite{diamandis2023efficient,angeris2022optimal}
with a network of $n$ assets and $m = (1/4)n^2$ markets between these assets.
Each market is a constant function market maker~\cite{angerisImprovedPriceOracles2020,angeris2023geometry}
that allows trades between a randomly selected pair of assets and has a strictly 
concave, increasing edge gain function. To interact with a market (\ie,
to use this edge), a trader must pay a fixed fee $q_0 \in \reals_+$. The trader's
goal is to maximize their utility $U(y)$ of the net trade $y$, given by
\[
    U(y) = c^Ty - \frac{\mu}{2}y^Ty,
\]
where $\mu > 0$ is some risk aversion parameter. The objective function is,
therefore,
\[
    c^Ty - \frac{\mu}{2}y^Ty + q_0(\ones^T\lambda)
\]
We solve the convex relaxation of the fixed fee problem~\eqref{eq:fixed-fee-relaxation}
using the open-source convex optimization solver Clarabel~\cite{clarabel}. All
code is available at
\begin{center}
    \texttt{\url{https://github.com/tjdiamandis/routing-theory-experiments}}.
\end{center}

\begin{table}[h]
    \begin{center}
        \begin{tabular}{@{} ll cccccc @{}}
            \toprule
            && \multicolumn{3}{c}{low fee ($q_0 = 0.01$)} & \multicolumn{3}{c@{}}{high fee $q_0 = 1.0$}\\
            \cmidrule(lr){3-5} \cmidrule(l){6-8} 
            $n$ & $m$ &$\mu = 0$ & $\mu = 10^{-4}$ & $\mu = 10^{-2}$ & $\mu = 0$ & $\mu = 10^{-4}$ & $\mu = 10^{-2}$  \\
            \midrule
            10&     25      &0      &0      &10     &0      &0      &10\\
            17&     72      &0      &0      &16     &0      &3      &16\\
            28&     196     &0      &0      &28     &0      &0      &26\\
            46&     529     &0      &0      &47     &0      &4      &46\\
            77&     1,482    &0      &0      &76     &0      &4      &75\\
            129&    4,160    &0      &3      &129    &0      &17     &126\\
            215&    11,556   &2      &11     &213    &0      &56     &212\\
            359&    32,220   &7      &45     &359    &1      &150    &355\\
            599&    89,700   &9      &536    &601    &3      &479    &592\\
            1,000&   250,000  &33     &1,143   &1,008   &5      &884    &992\\
            \bottomrule
        \end{tabular}
        \caption{Number of non-integral $\lambda_i$ for the fixed fee problem relaxation.}
        \label{tab:violations}
    \end{center}
\end{table}

\begin{table}[h]
    \begin{center}
        \begin{tabular}{@{} ll cccccc @{}}
            \toprule
            && \multicolumn{3}{c}{low fee ($q_0 = 0.01$)} & \multicolumn{3}{c@{}}{high fee $q_0 = 1.0$}\\
            \cmidrule(lr){3-5} \cmidrule(l){6-8} 
            $n$ & $m$ &$\mu = 0$ & $\mu = 10^{-4}$ & $\mu = 10^{-2}$ & $\mu = 0$ & $\mu = 10^{-4}$ & $\mu = 10^{-2}$  \\
            \midrule
            10      &25     &\texttt{8.20e-11}      &\texttt{1.10e-10}      &\texttt{9.11e-04}      &\texttt{1.78e-09}      &\texttt{4.26e-09}      &\texttt{9.32e-02}\\
            17      &72     &\texttt{3.44e-10}      &\texttt{7.65e-10}      &\texttt{8.83e-04}      &\texttt{1.92e-08}      &\texttt{6.04e-04}      &\texttt{9.37e-02}\\
            28      &196    &\texttt{1.12e-09}      &\texttt{6.43e-10}      &\texttt{6.99e-04}      &\texttt{7.21e-09}      &\texttt{4.77e-09}      &\texttt{7.01e-02}\\
            46      &529    &\texttt{2.28e-10}      &\texttt{1.15e-10}      &\texttt{7.61e-04}      &\texttt{1.23e-08}      &\texttt{2.60e-04}      &\texttt{8.37e-02}\\
            77      &1,482   &\texttt{2.77e-10}      &\texttt{1.46e-10}      &\texttt{8.42e-04}      &\texttt{5.85e-09}      &\texttt{5.77e-05}      &\texttt{9.21e-02}\\
            129     &4,160   &\texttt{4.27e-10}      &\texttt{1.57e-07}      &\texttt{9.33e-04}      &\texttt{2.66e-09}      &\texttt{1.05e-04}      &\texttt{1.03e-01}\\
            215     &11,556  &\texttt{3.10e-08}      &\texttt{4.06e-07}      &\texttt{8.27e-04}      &\texttt{6.73e-09}      &\texttt{1.87e-04}      &\texttt{9.20e-02}\\
            359     &32,220  &\texttt{6.70e-08}      &\texttt{6.72e-07}      &\texttt{8.39e-04}      &\texttt{1.16e-06}      &\texttt{2.60e-04}      &\texttt{9.25e-02}\\
            599     &89,700  &\texttt{2.47e-08}      &\texttt{5.36e-06}      &\texttt{8.91e-04}      &\texttt{4.22e-07}      &\texttt{4.85e-04}      &\texttt{9.82e-02}\\
            1,000    &250,000 &\texttt{3.73e-08}      &\texttt{7.05e-06}      &\texttt{9.20e-04}      &\texttt{7.65e-07}      &\texttt{5.48e-04}      &\texttt{1.02e-01}\\
            \bottomrule
        \end{tabular}
        \caption{Relative objective difference between the relaxation and the objective value with the rounded solution.}
        \label{tab:obj-diff}
    \end{center}
\end{table}

\paragraph{Results.}
We solve the fixed fee problem for logarithmically spaced range $n$ ranging from
$10$ to $1{,}000$ ($m = 25$ to $250{,}000$), for $\mu = 0$ (linear objective),
$10^{-4}$, and $10^{-2}$, and for both a `low fee' setting with $q_0 = 0.01$ while for
the high fee setting, we have $q_0 = 1.0$.
We record the number of non-integral $\lambda_i$ in the solution to the
relaxation in table~\ref{tab:violations} and the relative objective difference
between the relaxation and the objective value with the rounded solution,
using~\eqref{eq:feasible-approx}, in table~\ref{tab:obj-diff}. We see that when
the objective function is linear, the number of non-integral $\lambda_i$ is
negligible, and the objective value is very close to the optimal value as we
would expect. As the objective function becomes more concave, the number of
non-integral $\lambda_i$ increases, but remains small relative to $m$.
Additionally, we see that the objective value of the rounded solution is still
very close to that of the relaxation, suggesting that we can approximately solve
this NP-hard problem in practice. We leave further numerical investigation to
future work.

\section{Conclusion}
In this work, we have shown a number of theoretical properties of the convex
flow problem, all of which follow more or less directly from standard results in
convex geometry. We first showed that the traditional assumption of
nondecreasing objective terms can be replaced by a downward closedness condition
on the sets of allowable flows. Using this condition, we showed that there is a
natural calculus of flows. We then introduced the flow cone to derive a conic
formulation of the convex flow problem, which we showed is equivalent to the
original and has a number of interesting properties. We next examined the
case of fixed costs for the use of an edge. Via a Shapley--Folkman argument, we
showed that the relaxation of the fixed fee problem always has an `almost'
integral solution. Finally, we showed that this fixed fee problem can be solved
using the same algorithm for the convex flow problem, and solving the convex
relaxation often results in an integral solution that is close to optimal
in practice.

\paragraph{Future work.}
This work prompts a number of questions, any of which suggest an interesting
avenue for future research. First, why do we often find integral solutions to
the fixed fee problem? Is there a natural condition that guarantees this?
Second, what algorithm is best for solving the fixed fee problem? Should we use
the same algorithm as for the convex flow problem, with the modification
suggested in~\S\ref{sec:fixed-cost-dual}? Finally, this problem has a very
natural decomposition over the edges (equivalently, over the nodes, since we are
working with a bipartite graph). Can we exploit this decomposition to devise
efficient distributed, tatonnement-style algorithms? If so, what convergence
rates should we expect? These potentially asynchronous algorithms may be of
interest in decentralized applications, such as power grids or wireless networks
as discussed in~\cite[\S3]{diamandis2024convex}.

\section*{Acknowledgements}
Theo Diamandis is supported by the Department of Defense (DoD) through the 
National Defense Science \& Engineering Graduate (NDSEG) Fellowship Program.

\printbibliography

\appendix

\section{Downward closure of the flow cone}\label{app:downward-closure-flow-cone}
In this section, we show that the flow cone $K_i$ from~\eqref{eq:flow-cone} is
downward closed. 

\paragraph{Scaling property.} An important property of a set of allowable flows is
that, for any flow set $T_i$, if $x \in T_i$, then $\alpha x \in T_i$ for any
$0 \le \alpha \le 1$. The proof is nearly immediate by the convexity of $T_i$
and the fact that $0 \in T_i$ by noting that
\begin{equation}\label{eq:dilation}
    \alpha x = \alpha x + (1- \alpha)0 \in T_i.
\end{equation}

\paragraph{Downward closure.}
We will now prove that the flow cone $K_i$ is
downward closed. In other words, we will show that, for any $(x, \lambda) \in
K_i$ and $(x', \lambda')$ with $x' \le x$ and $\lambda' \le \lambda$, then
$(x', \lambda') \in K_i$. We will first show that $(x, \lambda) \in K_i$
implies that $(x, \lambda') \in K_i$. To see this, note that, if $(x, \lambda)
\in K_i$, then there exists a sequence $\{(x_k, \lambda_k)\}$ such that
$x_k/(-\lambda_k) \in T_i$ and $\lambda_k < 0$ for each $k$, while the
sequences converge in that $x_k \to x$ and $\lambda_k \to \lambda$ by the
definition~\eqref{eq:flow-cone}. If $\lambda = \lambda'$ then obviously $(x,
\lambda') \in K_i$, so assume that $\lambda' < \lambda$. In this case, there
exists some $k'$ such that, $\lambda_k \ge \lambda'$ for all $k \ge k'$. But,
since
\[
    \frac{x_k}{-\lambda_k} \in T_i,
\]
for all $k \ge k'$ by definition, then, since $0 \le \lambda_k/\lambda' \le 1$
we have that, using the scaling property~\eqref{eq:dilation}:
\[
    \frac{\lambda_k}{\lambda'}\frac{x_k}{-\lambda_k} \in T_i,
\]
so $x_k/(-\lambda') \in T_i$ for each $k \ge k'$. Since $x_k \to x$ and $T_i$
is closed, then $x/(-\lambda') \in T_i$, so, by definition $(x, \lambda') \in
K_i$, as required. The final question, if $x' \le x$ then $(x', \lambda') \in
K_i$ is easy: if $\lambda' < 0$ then $x/(-\lambda') \in T_i$ implies that
$x'/(-\lambda') \in T_i$ by the downward closure of $T_i$, so $(x', \lambda')
\in K_i$. If $\lambda' = 0$ then we know that $\lambda = 0$. Set $\delta = x' -
x$ and note that $\delta \le 0$ since $x' \le x$. This means that, for the
sequence $\lambda_k < 0$ with $\lambda_k \to 0$, we have
\[
    \frac{x_k + \delta}{-\lambda_k} \in T_i,
\]
since $x_k/(-\lambda_k) \in T_i$ and $(x_k + \delta)/(-\lambda_k) \le
x_k/(-\lambda_k)$ because $\delta \le 0$ and $T_i$ is downward closed. But,
since $x_k \to x$, then $x_k + \delta \to x + \delta = x'$. By the definition
of~\eqref{eq:flow-cone}, then we have that $(x + \delta, 0) \in K_i$, so $(x',
\lambda') = (x + \delta, 0) \in K_i$ as required.

\section{Downward closure equivalence using the dual}\label{sec:downward-closure-dual-proof}
Here, we provide a brief, intuitive proof of the equivalence between downward
closure and monotonicity using the dual problem.

\paragraph{Monotonicity implies downward closure.}
Assume that the objective functions $U$ and $\{V_i\}$ are nondecreasing. Then
the dual variables $\nu$ and $\{\eta_i\}$ must all be nonnegative for the dual
to be finite-valued. Consider the arbitrage problem for some $i$:
\[
    f_i(\eta_i) = \sup_{x_i \in T_i} \eta_i^Tx_i.
\]
The solution to the arbitrage problem can be broken into two cases. First, if
the set $T_i$ contains a positive ray, then the optimal value is unbounded, and
the primal problem is therefore infeasible. Second, if the optimal value is
finite, then the solution will lie on the boundary of $T_i$, as $T_i$ is closed
and convex. As a result, we can replace $T_i$ with its `downward extension',
\[
    \tilde T_i = T_i - \reals^n_+,
\]
without affecting the solution.

\paragraph{Downward closure implies monotonicity.}
Assume that the sets $\{T_i\}$ are downward closed. We use $x_i^\star$ to denote a
solution to the arbitrage problem:
\[
    x_i^\star = \argmax_{x \in T_i} \eta_i^Tx,
\]
for a fixed price vector $\eta_i$ when this solution is finite. Consider two
possible cases for $f_i$. First, if $f_i$ is unbounded, then the primal problem
is again infeasible. Second, if $f_i$ is is finite, then the set $T_i$ is
contained in the halfspace $\{x \mid \eta_i^Tx \le f_i(\eta_i)\}$, and $x_i^\star$
is on the boundary of $T_i$. Since $T_i$ is downward closed, we must have 
$\eta_i \ge 0$ in this case. As a result, we can replace the objective function 
with its monotonic concave envelope without affecting the solution.

\section{Other properties}
\subsection{Circulation problem}\label{app:circulation}
If the net flow utility simply constrains the net flow to be zero, \ie,
\[
   U(y) = -I_{\{0\}}(y).
\]
then we recover a generalized circulation problem. The dual problem becomes
\[
   \begin{aligned}
       & \text{minimize} && {\sum_{i=1}^m\left( \bar V_i(\eta_i - A_i^T\nu) + f_i(\eta_i)\right)}.
   \end{aligned}
\]
Taking the infimum of the objective over $\eta$ and introducing a new variables
$z_i \in \reals^{n_i}$ for $i = 1, \dots, m$, we can rewrite this problem as
\[
   \begin{aligned}
       & \text{minimize} && {\sum_{i=1}^m\tilde V_i(z_i)} \\
       & \text{subject to} && z_i = A_i^T\nu, \quad i=1, \dots, m\\
       &&& \nu \ge 0,
   \end{aligned}
\]
where we define
\[
   \tilde V_i(z_i) = \inf_{\eta_i} \left\{ \bar V_i(\eta_i - z_i) + f_i(\eta_i) \right\}.
\]
(This is very close to, but not quite, the infimal convolution of $\bar V_i$
and $f_i$.) Note that $\tilde V_i$ is a convex function, as convexity is
preserved under partial minimization~\cite[\S3.2.5]{cvxbook}. This problem has
a nice interpretation: we are finding the pin voltages on the $m$ components in
a passive electrical circuit~\cite[\S6]{boyd2007notes}. The optimality
conditions simplify to
\[
\begin{aligned}
   \nabla \tilde V_i(z_i) &= x_i, \quad i=1, \dots, m,\\
   A_i^T\nu &= z_i, \quad i=1, \dots, m, \\
   \sum_{i=1}^m A_ix_i &\ge 0, \quad \nu \ge 0.
\end{aligned}
\]
Viewing $z_i$ and $x_i$ as the voltage and current, respectively, at the
terminals of component $i$, and viewing $\nu$ as the voltages at every node in
the circuit, then the first equation can can be interpreted as the
voltage-current characteristic for component $i$, the second as the Kirchoff
voltage law, and the last as Kirchoff's current law, respectively.

\subsection{Cycle condition}
When all edges are between two nodes, the optimality conditions have a nice
interpretation in terms of cycle conditions, similar to the augmenting path
condition for max flow, given by Ford and Fulkerson. 

We define the vectors $\delta_1, \dots, \delta_m$ to be an arbitrage with 
respect to the flows $\{x_i\}$ if the following conditions hold:
\begin{enumerate}
    \item The vector $\delta_i \in T_i^\star(x_i)$ for all $i = 1, \dots, m$, where
    \[
        T_i^\star(x_i) = \{\delta \mid x_i + t\delta_i \in T_i ~\text{for some}~ t > 0\}.
    \]
    \item For some subgradient $\hat \nu \in U\left(\sum_{i=1}^m x_i \right)$,
    we have that
    \[
        \hat \nu^T\left(\sum_{i=1}^m A_i \delta_i \right) \le 0.
    \]
\end{enumerate}
For the standard max flow problem, this should give the classic augmenting path
condition. For more general problems, we can first decompose into a cycle basis
to easily check this condition. 

\end{document}